# Global Stability Analysis of the Age-Structured Chemostat With Substrate Dynamics

**Iasson Karafyllis[*], Dionysios Theodosis[*], and Miroslav Krstic[**]**

[*]Dept. of Mathematics, National Technical University of Athens,
Zografou Campus, 15780, Athens, Greece,
emails: iasonkar@central.ntua.gr , dtheodp@central.ntua.gr

[**]Dept. of Mechanical and Aerospace Eng., University of California, San Diego, La Jolla, CA 92093-0411, U.S.A., email: krstic@ucsd.edu

**Abstract**

In this paper we study the stability properties of the equilibrium point for an age-structured chemostat model with renewal boundary condition and coupled substrate dynamics under constant dilution rate. This is a complex infinite-dimensional feedback system. It has two feedback loops, both nonlinear. A positive static loop due to reproduction at the age-zero boundary of the PDE, counteracted and dominated by a negative dynamic loop with the substrate dynamics. The derivation of explicit sufficient conditions that guarantee global stability estimates is carried out by using an appropriate Lyapunov functional. The constructed Lyapunov functional guarantees global exponential decay estimates and uniform global asymptotic stability with respect to a measure related to the Lyapunov functional. From a biological perspective, stability arises because reproduction is constrained by substrate availability, while dilution, mortality, and substrate depletion suppress transient increases in biomass before age-structure effects can amplify them. The obtained results are applied to a chemostat model from the literature, where the derived stability condition is compared with existing results that are based on (necessarily local) linearization methods.

## 1. Introduction

Chemostat models constitute a fundamental framework for describing microbial growth in a bioreactor with continuous nutrient inflow and outflow. Classical chemostat models use a single variable to represent the microorganisms leading to systems of ordinary differential equations (see for instance [1], [3], [25], [29]). However, in many biological systems the evolution of the population of microorganisms depends strongly on their age. Incorporating age structure leads to infinite-dimensional models described by first-order hyperbolic partial differential equations with renewal boundary conditions like the well-known McKendrick-von Foerster equation (see [7], [23]). This additional structure substantially alters the qualitative behavior of the system. In contrast to classical ODE formulations with constant yield coefficients, age-structured models may exhibit oscillatory behavior that can reproduce experimentally observed data (see [26], [32], [33]).

When the age distribution is coupled with substrate dynamics, the mathematical structure becomes considerably more involved. The resulting system consists of a transport Partial Differential Equation (PDE) for the population age density and a nonlinear ordinary differential equation for the substrate, coupled with the PDE through nonlocal terms (see for instance [18], [32], [33]). This nonlinear and nonlocal coupling introduces substantial analytical difficulties, especially in the study of global qualitative properties such as boundedness, attractivity, and asymptotic stability.

For chemostat models described by finite-dimensional systems, equilibrium analysis and stability properties are well understood (see for instance, [1], [5], [27], [29]). In contrast, for age-structured chemostats, existing results are sparse and deal only with local stability properties (obtained via linearization techniques), or deal with settings in which the substrate equation is absent (see for instance [11], [12], [13], [10], [24], [25], [32], [33]). Establishing global stability for the age-structured chemostat with substrate dynamics therefore remains substantially more challenging, due to the nonlinear and nonlocal interactions and the positivity constraints inherent to population dynamics.

The objective of this work is to provide sufficient conditions for uniform global asymptotic stability of the equilibrium point for an age-structured chemostat model. We provide global $KL$ stability estimates for an age-structured chemostat model coupled with substrate dynamics under constant dilution rate.

Though no feedback law is being designed in the paper, the paper conducts nonlinear infinite-dimensional feedback analysis. The model has two nonlinear feedback loops. One is a positive static feedback loop due to reproduction at the birth boundary condition of the PDE, while the other loop is dynamic and negative through the substrate dynamics. The main analytical challenge lies in quantifying, using tools that extend far beyond the conventional small-gain theorem concepts, the interaction between these two loops and establishing conditions under which the negative feedback dominates.

To address this challenge, our approach leverages the recent results obtained in [18], that established global existence of solutions. Under suitable structural assumptions linking the birth and substrate consumption rates, we first show the existence of a trapping (absorbing) region that is reached by all solutions in finite time. We then construct a Lyapunov functional and derive explicit sufficient conditions that guarantee a global $KL$ stability estimate with respect to a specific measure. The construction of the Lyapunov functional combines deviations of the age profile with logarithmic transformations of the normalized substrate and reproduction rates, allowing us to deal with the nonlinear and nonlocal terms of the model. The result establishes global exponential decay estimates and convergence of the age distribution together with the substrate concentration. The provided sufficient conditions for global asymptotic stability are explicit and can be readily verified from the model data, making them applicable in a wide range of cases. In intuitive and biological terms, we prove that the equilibrium is globally asymptotically stable, under suitable conditions, due to reproduction being constrained by substrate availability, and due to the combined effects of dilution, mortality, and substrate depletion providing sufficient negative feedback to damp any temporary increase in biomass faster than age-cohort effects can build up self-sustaining oscillations. Finally, as an illustration, we apply the derived stability theorem to a specific age-structured chemostat model studied in [32], for which local attractivity was established via linearization, and compare the local criteria given in [32] with the global stability condition obtained in our work.

The paper is organized as follows. Section 2 introduces the model and presents its basic properties. Section 3 provides the main results, including the existence of a trapping region and the global $KL$ stability estimates. Section 4 provides an illustrative example and comparison with existing results. Section 5 contains the proofs of the main results. Finally, some concluding remarks are given in Section 6.

**Notation** Throughout this paper, we adopt the following notation.

* $\mathbb{R}_+ := [0,+\infty)$. For a vector $x \in \mathbb{R}^n$, $|x|$ denotes its Euclidean norm. We use the notation $x^+$ for the positive part of the real number $x \in \mathbb{R}$, i.e., $x^+ = \max(x,0)$.

* Let $A \subseteq \mathbb{R}^n$ be an open set and let $B \subseteq \mathbb{R}^n$ be a set that satisfies $A \subseteq B \subseteq cl(A)$, where $cl(A)$ is the closure of $A$. By $C^0(B;\Omega)$, we denote the class of continuous functions on $B$, which take values in $\Omega \subseteq \mathbb{R}^m$. By $C^k(B;\Omega)$, where $k \ge 1$ is an integer, we denote the class of functions on $B \subseteq \mathbb{R}^n$, which take values in $\Omega \subseteq \mathbb{R}^m$ and have continuous derivatives of order $k$. In other words, the functions of class $C^k(B;\Omega)$ are the functions which have continuous derivatives of order $k$ in $A = \mathrm{int}(B)$ that can be continued continuously to all points in $\partial A \cap B$. When $\Omega = \mathbb{R}$ then we write $C^0(B)$ or $C^k(B)$.

* Let $A \subseteq \mathbb{R}^n$ be an open set and let $\Omega \subseteq \mathbb{R}^m$ be a non-empty set. By $L^p(A;\Omega)$ with $p \ge 1$ we denote the equivalence class of measurable functions $f : A \to \Omega$ for which $\|f\|_p = \left( \int_A |f(x)|^p \, dx \right)^{1/p} < +\infty$. By $L^\infty(A;\Omega)$ we denote the equivalence class of

measurable functions $f: A \to \Omega$ for which $\|f\|_\infty = \sup_{x\in A}\left(|f(x)|\right) < +\infty$ where $\sup_{x\in A}\left(|f(x)|\right)$ is the essential supremum. When $\Omega = \mathbb{R}^m$ we simply write $L^p(A)$. When $B \subseteq \mathbb{R}^n$ is not open but has non-empty interior, $L^p(B;\Omega)$ and $L^\infty(B;\Omega)$ mean $L^p(A;\Omega)$ and $L^\infty(A;\Omega)$, respectively, with $A = \mathrm{int}(B)$. By $L^\infty_{loc}(\mathbb{R}_+;\Omega)$ we denote the equivalence class of measurable functions $f:\mathbb{R}_+ \to \Omega$ with $f \in L^\infty([0,T];\Omega)$ for every $T>0$.

* By $K$ we denote the class of increasing $C^0$ functions $a:\mathbb{R}_+ \to \mathbb{R}_+$ with $a(0)=0$. By $K_\infty$ we denote the class of increasing $C^0$ functions $a:\mathbb{R}_+ \to \mathbb{R}_+$ with $a(0)=0$ and $\lim_{s\to+\infty} a(s) = +\infty$. By $KL$ we denote the set of all continuous functions $\sigma:\mathbb{R}_+ \times \mathbb{R}_+ \to \mathbb{R}_+$ with the following properties: (i) for each $t \ge 0$ the mapping $\sigma(\cdot,t)$ is of class $K$; (ii) for each $s \ge 0$, the mapping $\sigma(s,\cdot)$ is non-increasing with $\lim_{t\to+\infty} \sigma(s,t) = 0$.

* Let $I \subseteq \mathbb{R}_+$ and let $f: I \times \mathbb{R}_+ \to \mathbb{R}_+$ be a given function. We use the notation $f[t]$ to denote the profile at certain $t \in I$, i.e., $(f[t])(a) = f(t,a)$ for all $a \ge 0$.

* Let $k \in L^\infty(\mathbb{R}_+)$ and $f \in L^1(\mathbb{R}_+)$ be given. We define $\langle k, f\rangle := \int_0^{+\infty} k(a) f(a)\, da$.

## 2. The Age-Structured Chemostat Model

Consider the age-structured chemostat model

$$\frac{\partial f}{\partial t}(t,a) + \frac{\partial f}{\partial a}(t,a) = -\left(\beta(a) + D(t)\right) f(t,a) \tag{2.1}$$

$$f(t,0) = \mu\left(S(t)\right)\langle k, f[t]\rangle \tag{2.2}$$

$$\dot{S}(t) = D(t)\left(S_{in} - S(t)\right) - \mu\left(S(t)\right)\langle q, f[t]\rangle \tag{2.3}$$

where $f(t,a) > 0$ is the distribution function of the microbial population in the chemostat at time $t \ge 0$ and age $a \ge 0$, $S(t) > 0$ is the limiting substrate concentration, $S_{in} > 0$ is the inlet concentration of the substrate, $D(t) > 0$ is the dilution rate, $\mu(S)$ is the specific growth rate function, $\beta(a)$ is the mortality rate and $k(a), q(a)$ are functions that determine the birth of new cells and the substrate consumption of the microbial population, respectively. All functions $\mu, \beta, k, q: \mathbb{R}_+ \to \mathbb{R}_+$ are assumed to be bounded, $C^0(\mathbb{R}_+)$ functions with $\mu \in C^1(\mathbb{R}_+)$, $\mu(0)=0$, $\mu(S) > 0$ for $S > 0$ and $\int_0^{+\infty} k(a)\, da > 0$, $\int_0^{+\infty} q(a)\, da > 0$.

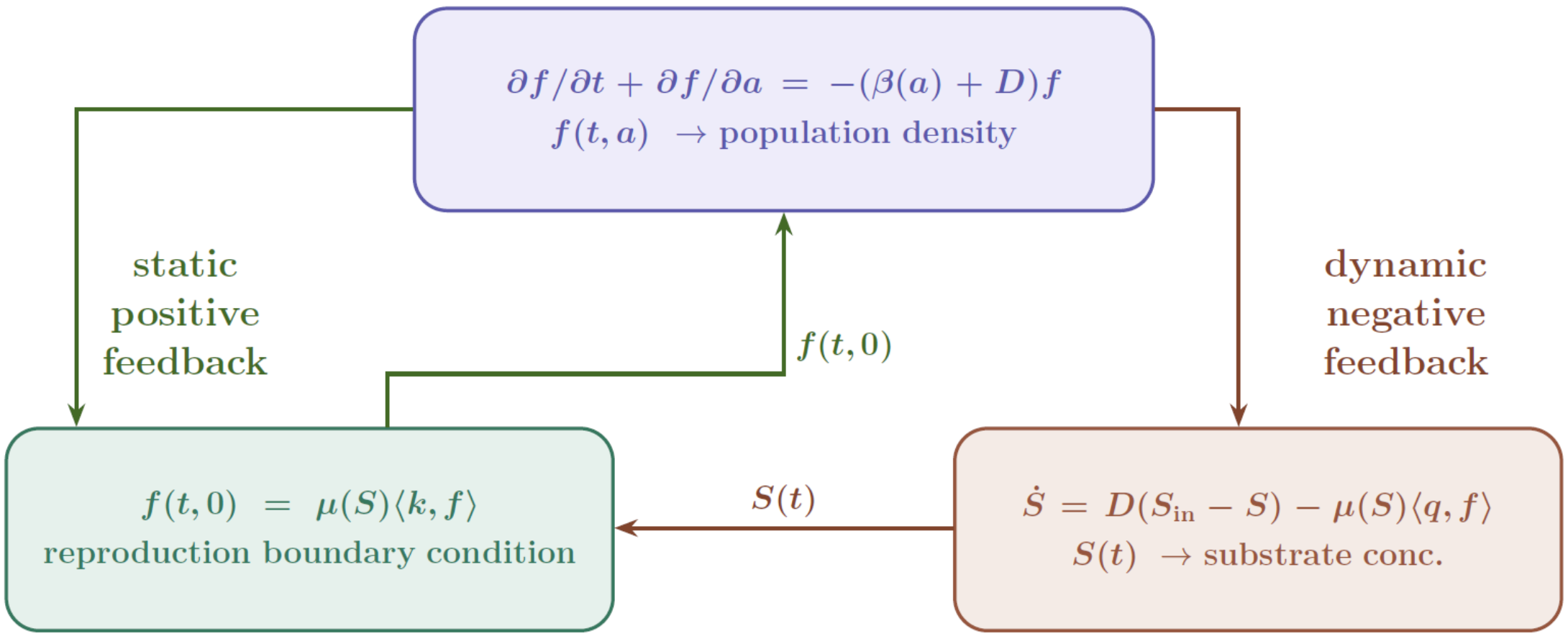

**Figure 1:** Feedback structure of the age-structured chemostat model (2.1), (2.2), (2.3).

The interconnection structure of system (2.1), (2.2), (2.3) is illustrated in Figure 1. From a control systems perspective, the model can be viewed as the interconnection of two feedback loops. The first loop is a positive feedback associated with the boundary condition (2.2), where the reproduction rate $\langle k, f \rangle$ drives the inflow of newborns $f(t,0)$ and, hence, also the subsequent age distribution $f(t,a)$. The second loop is a dynamic negative feedback through the substrate equation (2.3), where the consumption $\langle q, f \rangle$ reduces the substrate concentration $S(t)$, which in turn limits reproduction via the growth function $\mu(S)$. This feedback interaction plays a central role in the stability analysis.

Define:

$$X = \left\{ (f,S) \in L^1(\mathbb{R}_+) \times (0, S_{in}) : \begin{array}{c} f \in C^0(\mathbb{R}_+;(0,+\infty)),\ \lim\limits_{a \to +\infty} \left( f(a) \right) = 0 \\ f(0) = \mu(S) \int\limits_0^{+\infty} k(a) f(a) da \end{array} \right\} \quad (2.4)$$

We consider $X$ to be a metric space with metric given by the formula for all $(f,S), (\bar{f}, \bar{S}) \in X$:

$$d\left( (f,S), (\bar{f}, \bar{S}) \right) = \| f - \bar{f} \|_1 + |S - \bar{S}| \quad (2.5)$$

We next provide the notion of solution that is appropriate for system (2.1), (2.2), (2.3).

**Definition 1:** *Let $D \in L^\infty_{loc}(\mathbb{R}_+;\mathbb{R}_+)$, $(f_0, S_0) \in X$ and $T > 0$ be given. We say that a continuous mapping $(f,S):[0,T] \to X$ is a weak solution on $[0,T]$ with input $D$ of the initial-boundary value problem (2.1), (2.2), (2.3) with initial condition*

$$f[0] = f_0\ ,\ S(0) = S_0 \quad (2.6)$$

*if the following properties are valid:*
*i) (2.6) holds,*

*ii)* $f \in C^0\left([0,T]\times\mathbb{R}_+;(0,+\infty)\right)$ *and* $S:[0,T]\to(0,S_{in})$ *is absolutely continuous,*

*iii) (2.2) holds for all* $t\in[0,T]$ *and (2.3) holds for* $t\in[0,T]$ *a.e.*

*iv) the following equation holds for all* $\phi \in C^1\left([0,T]\times\mathbb{R}_+\right)\cap L^\infty\left([0,T]\times\mathbb{R}_+\right)$ *with* $\left(\frac{\partial\phi}{\partial a}+\frac{\partial\phi}{\partial t}\right)\in L^\infty\left([0,T]\times\mathbb{R}_+\right)$ *and* $t\in[0,T]$*:*

$$\begin{aligned}&\int_0^{+\infty} f_0(a)\phi(0,a)da+\int_0^t f(s,0)\phi(s,0)ds=\int_0^{+\infty} f(t,a)\phi(t,a)da\\&+\int_0^t\int_0^{+\infty}\left(\left(\beta(a)+D(s)\right)\phi(s,a)-\frac{\partial\phi}{\partial a}(s,a)-\frac{\partial\phi}{\partial s}(s,a)\right)f(s,a)dads\end{aligned}\tag{2.7}$$

*We say that the initial-boundary value problem (2.1), (2.2), (2.3), (2.6) admits a global weak solution with input* $D$ *if for every* $T>0$ *there exists a weak solution on* $[0,T]$ *with input* $D$ *of the initial-boundary value problem (2.1), (2.2), (2.3), (2.6).*

Next, we highlight a list of properties for the age-structured chemostat model (2.1), (2.2), (2.3) and its solutions.

**Properties**:

1) Theorem 1 in [18] shows that for any $(f_0,S_0)\in X$, $D\in L^\infty_{loc}\left(\mathbb{R}_+;\mathbb{R}_+\right)$ the weak solution $\left(f[t],S(t)\right)$ with input $D$ of the initial-boundary value problem (2.1), (2.2), (2.3), (2.6) is unique and is defined for all $t\geq 0$. Notice that for constant dilution rate $D(t)\equiv D>0$, $S:\mathbb{R}_+\to(0,S_{in})$ is continuously differentiable.

2) Definition (2.4) of the state space $X$ (the facts that $f\in C^0\left(\mathbb{R}_+;(0,+\infty)\right)$ and $\lim_{a\to+\infty}\left(f(a)\right)=0$), implies that $X\subseteq L^\infty\left(\mathbb{R}_+\right)\times\mathbb{R}$. Consequently, (since $\|f\|_p\leq\|f\|_1^{1/p}\|f\|_\infty^{1-1/p}$ for every $p\in[1,+\infty]$ and every $f\in L^\infty\left(\mathbb{R}_+\right)\cap L^1\left(\mathbb{R}_+\right)$) it follows that $X\subseteq L^p\left(\mathbb{R}_+\right)\times\mathbb{R}$ for all $p\in[1,+\infty]$. Thus, for any $(f_0,S_0)\in X$, $D\in L^\infty_{loc}\left(\mathbb{R}_+;\mathbb{R}_+\right)$ the weak solution $\left(f[t],S(t)\right)$ with input $D$ of the initial-boundary value problem (2.1), (2.2), (2.3), (2.6) satisfies for all $t\geq 0$ and $p\in[1,+\infty]$:

$$f[t]\in L^\infty\left(\mathbb{R}_+\right)\tag{2.8}$$

$$\|f[t]\|_p\leq\|f[t]\|_1^{1/p}\|f[t]\|_\infty^{1-1/p}\tag{2.9}$$

3) Let $D(t)\equiv D>0$ (a constant) be given and assume that the growth rate function $\mu(S)$ is increasing. An equilibrium point for (2.1), (2.2), (2.3) that corresponds to the (constant) dilution rate $D$ consists of a point $S^* \in (0, S_{in})$ that satisfies the so-called Lotka-Sharpe condition

$$\mu\left(S^*\right)\langle k, r\rangle = 1 \tag{2.10}$$

and has the equilibrium distribution function given for $a \geq 0$ by the formula

$$f^*(a) = \frac{D\left(S_{in} - S^*\right)}{\mu\left(S^*\right)\langle q, r\rangle} r(a) \tag{2.11}$$

where

$$r(a) = \exp\left( -Da - \int_0^a \beta(s)ds \right) \tag{2.12}$$

4) If we assume that $k \in C^1(\mathbb{R}_+) \cap L^\infty(\mathbb{R}_+)$ with $k' \in L^\infty(\mathbb{R}_+)$ then by taking $\phi(t,a) = k(a)$ we obtain from (2.7) for $D(t) \equiv D > 0$:

$$\begin{aligned} \langle k, f_0\rangle + k(0)\int_0^t f(s,0)ds &= \langle k, f[t]\rangle + \int_0^t \int_0^{+\infty} \left( \left(\beta(a) + D\right)k(a) - k'(a)\right) f(s,a)dads \\ &= \langle k, f[t]\rangle + \int_0^t \langle \beta k, f[s]\rangle ds + D\int_0^t \langle k, f[s]\rangle ds - \int_0^t \langle k', f[s]\rangle ds \end{aligned} \tag{2.13}$$

Since the map $(f,S): \mathbb{R}_+ \to X$ is continuous, it follows from (2.5) that the map $f: \mathbb{R}_+ \to L^1(\mathbb{R}_+)$ is continuous. It follows from (2.13) and the facts that $t \mapsto f(t,0) \in \mathbb{R}$ is continuous and $\beta, k, k' \in L^\infty(\mathbb{R}_+)$, that the map $t \mapsto \langle k, f[t]\rangle$ is continuously differentiable and in conjunction with (2.2) satisfies the following equation for all $t \geq 0$:

$$\frac{d}{dt}\left(\langle k, f[t]\rangle\right) = \left(k(0)\mu\left(S(t)\right) - D\right)\langle k, f[t]\rangle - \langle \beta k, f[t]\rangle + \langle k', f[t]\rangle \tag{2.14}$$

Similarly, if we assume that $q \in C^1(\mathbb{R}_+) \cap L^\infty(\mathbb{R}_+)$ with $q' \in L^\infty(\mathbb{R}_+)$ then by taking $\phi(t,a) = q(a)$ we obtain from (2.2), (2.7) for $D(t) \equiv D > 0$ that the map $t \mapsto \langle q, f[t]\rangle$ is continuously differentiable and satisfies the following equation for all $t \geq 0$:

$$\frac{d}{dt}\left(\langle q, f[t]\rangle\right) = q(0)\mu\left(S(t)\right)\langle k, f[t]\rangle - D\langle q, f[t]\rangle - \langle \beta q, f[t]\rangle + \langle q', f[t]\rangle \tag{2.15}$$

Finally, by taking $\phi(t,a)\equiv 1$ and recalling that $f(t,a)>0$ for all $t,a\geq 0$, we obtain from (2.2), (2.7) for $D(t)\equiv D>0$ that the map $t\mapsto \|f[t]\|_1$ is continuously differentiable and satisfies the following equation for all $t\geq 0$:

$$\frac{d}{dt}\left(\|f[t]\|_1\right)=\mu\left(S(t)\right)\langle k,f[t]\rangle - D\|f[t]\|_1-\langle \beta,f[t]\rangle \tag{2.16}$$

The following lemma provides some rough estimates for the weak solution $\left(f[t],S(t)\right)$ with constant input $D$ of the initial-boundary value problem (2.1), (2.2), (2.3), (2.6). Its proof is provided in Section 5.

**Lemma 1:** *Let $\left(f^*,S^*\right)\in X$ be an equilibrium point for (2.1), (2.2), (2.3) that corresponds to the constant input $D>0$. Then there exist constants $\bar{R},\bar{C}>0$ such that for every $\left(f_0,S_0\right)\in X$ the weak solution $\left(f[t],S(t)\right)$ with constant input $D$ of the initial-boundary value problem (2.1), (2.2), (2.3), (2.6) satisfies the following estimates for all $t\geq 0$:*

$$d\left(\left(f[t],S(t)\right),\left(f^*,S^*\right)\right)\leq \exp\left(\bar{R}t\right)d\left(\left(f_0,S_0\right),\left(f^*,S^*\right)\right) \tag{2.17}$$

$$\left\|f[t]-f^*\right\|_\infty \leq \max\left(\left\|f_0-f^*\right\|_\infty,\bar{C}\exp\left(\bar{R}t\right)d\left(\left(f_0,S_0\right),\left(f^*,S^*\right)\right)\right) \tag{2.18}$$

$$S(t)\leq S_{in}-\left(S_{in}-S_0\right)\exp\left(-Dt\right) \tag{2.19}$$

## 3. Main Results

In this section we establish the main qualitative properties of the solutions of system (2.1)–(2.3) with constant dilution $D(t)\equiv D>0$. We first show the existence of an absorbing (trapping) region, that is reached by all solutions in finite time. We then investigate the stability properties of the equilibrium introduced in Section 2 and prove convergence of solutions toward this equilibrium.

We use the following assumptions:

**(A)** *There exists a constant $R>0$ such that*

$$k(a)\leq Rq(a) \text{ for all } a\geq 0 \tag{3.1}$$

**(B)** *The function $\mu\in C^1\left(\mathbb{R}_+\right)$ is increasing, $k,q:\mathbb{R}_+\to\mathbb{R}_+$ are bounded, $C^1\left(\mathbb{R}_+\right)$ functions with $k',q'\in L^\infty\left(\mathbb{R}_+\right)$, and there exist constants $b,\gamma\in\mathbb{R}$, $\theta>0$, $\alpha,\delta\geq 0$ and continuous functions $h,p:\mathbb{R}_+\to\mathbb{R}_+$ such that*

$$
\begin{aligned}
&q'(a)-\beta(a)q(a)-\gamma Dq(a)-\delta\theta k(a)=h(a)\\
&k'(a)-\beta(a)k(a)-\frac{\alpha}{\theta}q(a)-bDk(a)=p(a)
\end{aligned}
\tag{3.2}
$$

Assumption (A) requires that the ratio of the birth modulus $k$ to the substrate consumption rate $q$ is bounded from above. Such an assumption is justified from a biological point of view, since only a part of the substrate consumption is utilized for reproduction while the major part of substrate consumption is used for endurance. Assumption (B) imposes boundedness on the derivatives of the birth and substrate consumption rates which reflects the fact that reproduction and substrate consumption vary smoothly with age and do not change abruptly as microorganisms age. Moreover, assumption (B) introduces certain relations between $k(a)$ and $q(a)$ coupling their evolution with the mortality $\beta(a)$ and the dilution $D$, where $h(a), p(a)$ account for additional age-dependent factors that may affect the substrate consumption and the reproduction. Assumptions similar to assumption (B) have been used in the literature (see [15], [32]).

The following lemma shows the existence of a trapping region. Its proof is provided in Section 5.

**Lemma 2:** *Suppose that assumptions (A), (B) hold. Then for every $F > RS_{in}$, where $R > 0$ is the constant involved in assumption (A), there exists an increasing continuous function $T:\mathbb{R}_+ \to \mathbb{R}_+$ such that for every $(f_0, S_0) \in X$ the weak solution $(f[t], S(t))$ of the initial-boundary value problem (2.1), (2.2), (2.3), (2.6) satisfies for all $t \geq T\left(RS_0 + \|f_0\|_1\right)$*

$$RS(t)+\|f[t]\|_1 \leq F \text{ and } S(t) \geq \underline{S} \tag{3.3}$$

*where*

$$\underline{S} = \frac{DS_{in}}{2(D+L_\mu F\|q\|_\infty)} \tag{3.4}$$

*and*

$$L_\mu = \max\left\{ \mu'(S) : S \in [0, S_{in}] \right\} \tag{3.5}$$

*Moreover, if $RS_0 + \|f_0\|_1 \leq F$ and $S_0 \geq \underline{S}$ then (3.3) holds for all $t \geq 0$.*

Let $F > RS_{in}$ be given and let $\underline{S} > 0$ be defined by (3.4). We define the set:

$$\Omega := \left\{ (f,S) \in X : RS + \|f\|_1 \leq F, S \geq \underline{S} \right\} \tag{3.6}$$

It is clear that Lemma 2 guarantees the set $\Omega$ defined by (3.6) is a trapping region of the age-structured chemostat model (2.1), (2.2), (2.3).

Next, we introduce the scalar variables

$$\zeta := \frac{\langle q, f \rangle}{\langle q, f^* \rangle}, \xi := \frac{\langle k, f \rangle}{\langle k, f^* \rangle}, \tag{3.7}$$

which represent the normalized substrate consumption and reproduction rates, respectively. It should be noticed that definition (2.4) and the fact that $k, q : \mathbb{R}_+ \to \mathbb{R}_+$ are continuous functions with $\int_0^{+\infty} k(a)da > 0$, $\int_0^{+\infty} q(a)da > 0$ guarantees that $\langle q, f \rangle > 0$ and $\langle k, f \rangle > 0$ for all $(f, S) \in X$. Therefore, definition (3.7) implies that $\xi, \zeta > 0$ for all $(f, S) \in X$.

To study the stability of the equilibrium $(f^*, S^*)$, we introduce the following normalized variables

$$v(a) := \frac{f(a)}{f^*(0)\zeta}, \; \chi(a) := v(a) - r(a) \tag{3.8}$$

where $r(a)$ is given by (2.12) and $\zeta$ given by (3.7). To capture the deviations of the substrate concentration from the equilibrium we also define

$$\varphi := \ln\left(\frac{\xi}{\zeta}\right), \; w := \ln\left(\frac{S_{in} - S^*}{S_{in} - S}\zeta\right) \tag{3.9}$$

Let also

$$g(S) = \frac{\mu(S)}{\mu(S^*)}, \kappa_1 = \frac{q(0)}{\langle q, r \rangle}, \kappa_2 = \frac{k(0)}{\langle k, r \rangle} \tag{3.10}$$

Under assumption (B) equations (2.14), (2.15) are valid. Combining (2.2), (2.3), (2.14), (2.15) (3.7), (3.8), (3.9), and (3.10), we obtain the equations

$$\chi(0) = g(S)\left(\exp(\varphi) - 1\right) + g(S) - 1 \tag{3.11}$$

$$\begin{aligned} \dot{\varphi} = {} & \left(\alpha + \theta \frac{\langle p, v \rangle}{\langle q, r \rangle}\right)\left(\exp(-\varphi) - 1\right) - \left(\kappa_1 g(S) + \delta\right)\left(\exp(\varphi) - 1\right) \\ & + (\kappa_2 - \kappa_1)\left(g(S) - 1\right) + \frac{\langle \theta p - h, \chi \rangle}{\langle q, r \rangle} \end{aligned} \tag{3.12}$$

$$\dot{w} = \left(\kappa_1 g(S) + \delta\right)\left(\exp(\varphi) - 1\right) - Dg(S)\left(\exp(w) - 1\right) + (\kappa_1 - D)\left(g(S) - 1\right) + \frac{\langle h, \chi \rangle}{\langle q, r \rangle} \tag{3.13}$$

$$\dot{S} = -D(S_{in} - S)\left(g(S) - 1\right) - D(S_{in} - S)g(S)\left(\exp(w) - 1\right) \tag{3.14}$$

For $(f,S)\in X$, and $\zeta,\xi,\chi,\varphi,w$ defined by (3.7), (3.8), and (3.9), respectively, we introduce the following Lyapunov functional

$$V(f,S):=\exp(\varphi)-\varphi-1+\Gamma\left(\exp(w)-w-1\right)+Q(S)+\frac{B}{2}\int_{0}^{+\infty}\rho^{2}(a)\chi^{2}(a)da \qquad (3.15)$$

where

$$\rho(a)=\exp(-\sigma a) \qquad (3.16)$$

and $Q:(0,S_{in})\to\mathbb{R}_{+}$ is a continuously differentiable function that satisfies

$$Q'(S):=M\frac{\left(g(S)-1\right)}{(S_{in}-S)g(S)} \qquad (3.17)$$

while $B,\Gamma,M>0,\sigma\geq 0$ are constants which are to be selected. The functional $V:X\to\mathbb{R}_{+}$ defined by (3.15) is well defined for every $(f,S)\in X$.

The following theorem is the main result of the paper. Theorem 1 provides sufficient conditions for a global $KL$ stability estimate with respect to the functional $\Psi:X\to\mathbb{R}_{+}$ defined by

$$\Psi(f,S):=\left|S-S^{*}\right|+\left\|f-f^{*}\right\|_{\infty}+\left\|f-f^{*}\right\|_{1}+V(f,S) \qquad (3.18)$$

It should be noted that the functional $\Psi(f,S)$ defined by (3.18) has many features of a measure or more accurately, the analogue of the notion of a measure (or size function) used for finite-dimensional systems (see for instance [30], [31]) defined on open sets. More specifically, the functional $\Psi(f,S)$ defined by (3.18) satisfies $\Psi\left(f^{*},S^{*}\right)=0$ and $\Psi(f,S)>0$ for all $(f,S)\in X\setminus\left\{\left(f^{*},S^{*}\right)\right\}$ and (as shown in Section 5 below) $\Psi(f,S)\to+\infty$ as $S\to 0^{+}$ or $S\to S_{in}^{-}$ or $\|f\|_{1}\to 0^{+}$ or $\left\|f-f^{*}\right\|_{p}\to+\infty$ for some $p\in[1,+\infty]$ (recall (2.9)). Therefore, the functional $\Psi(f,S)$ "blows up" when the state $(f,S)\in X$ approaches the "boundary" of the state space $X$.

**Theorem 1:** *Let $F>RS_{in}$ be given, let $\underline{S}>0$ be defined by (3.4) and let $\Omega$ be the trapping region defined by (3.6). Assume that assumptions (A) and (B) hold with*

$$\theta:=\mu\left(S^{*}\right)\langle q,r\rangle=\frac{\langle q,r\rangle}{\langle k,r\rangle} \qquad (3.19)$$

*Let* $\sigma \geq 0$ *be a constant and define* $\rho(a)$ *by (3.16). Let* $L := \inf_{a \geq 0} \left( \beta(a) \right)$ *and suppose that there exist constants,* $R_1, R_2, R_3 > 0$*,* $\omega > 0$*,* $\varepsilon > 0$*,* $\lambda \in [0,1]$*,* $B, \Gamma, M > 0$ *such that*

$$2\langle q,r\rangle \alpha \geq \omega\lambda \tag{3.20}$$

$$\kappa_1 g(\underline{S}) + \delta \geq \frac{\lambda \left\| \rho^{-1}(\theta p - h) \right\|_2^2}{2\omega B \langle q,r\rangle} \tag{3.21}$$

$$L + \gamma D + \sigma > \frac{\left\| \rho^{-1} h \right\|_2}{\langle q,r\rangle} \left( \left\| \rho r \right\|_2 + \frac{\Gamma}{2BR_2} \right) + \frac{(1-\lambda) \left\| \rho^{-1}(\theta p - h) \right\|_2}{2BR_1 \langle q,r\rangle} + \frac{\left\| \rho r \right\|_2^2}{2\varepsilon} \tag{3.22}$$

$$1 > 2B\varepsilon\delta \tag{3.23}$$

*and the following matrix is positive definite:*

$$P = \begin{bmatrix} A & -\Gamma\kappa_1/2 & (\kappa_1 - \kappa_2)/2 \\ -\Gamma\kappa_1/2 & \Gamma\left( D - \frac{R_2 \left\| \rho^{-1} h \right\|_2}{2g(\underline{S})\langle q,r\rangle} - \frac{\delta}{2R_3 g(\underline{S})} \right) & (\Gamma D + DM - \Gamma\kappa_1)/2 \\ (\kappa_1 - \kappa_2)/2 & (\Gamma D + DM - \Gamma\kappa_1)/2 & DM - Bg(S_{in})\left(1 + \varepsilon\kappa_1^2\right) \end{bmatrix} \tag{3.24}$$

*where*

$$A = \kappa_1 (1 - 2B\varepsilon\delta) + \frac{\delta(1 - B\varepsilon\delta)}{g(S_{in})} - B\left(1 + \varepsilon\kappa_1^2\right) g(S_{in})$$
$$- \frac{\Gamma\delta R_3}{2g(\underline{S})} - \frac{(1-\lambda) R_1}{2\langle q,r\rangle g(\underline{S})} \left\| \rho^{-1}(\theta p - h) \right\|_2$$

*Then, there exist a function* $C \in K_\infty$ *and a non-increasing function* $\bar{\sigma} : \mathbb{R}_+ \to (0, +\infty)$ *such that for every* $(f_0, S_0) \in X$ *the weak solution of the initial-boundary value problem (2.1), (2.2), (2.3), (2.6) satisfies the following estimate for all* $t \geq 0$*:*

$$\Psi(f[t], S(t)) \leq \exp\left( -\bar{\sigma}\left( C\left( \Psi(f_0, S_0) \right) \right) t \right) C\left( \Psi(f_0, S_0) \right) \tag{3.25}$$

The proof of Theorem 1 is provided in Section 5.

**Remarks on Theorem 1:**

**(a)** The functional $\Psi(f, S)$ can be interpreted as a composite measure of the deviation from equilibrium. In particular, the deviation of $f$ from $f^*$ is captured directly in the $L^1$ and $L^\infty$ norms, through a weighted $L^2$- type term measuring the deviation of a normalized age profile, and through quantities related to the birth and consumption rates. More precisely, $\Psi(f,S)$ contains the direct terms $|S-S^*|$, $\|f-f^*\|_1$, and $\|f-f^*\|_\infty$, together with the Lyapunov functional $V(f,S)$. The latter includes a weighted quadratic term measuring the deviation of the normalized age profile $v(a)$ from the normalized equilibrium profile $r(a)$, as well as the logarithmic terms $e^\varphi-\varphi-1$ and $e^w-w-1$, where $\varphi=\ln(\xi/\zeta)$ measures the deviation between normalized birth and normalized consumption rates, and $w$ compares the substrate deficit ( $S_{in}-S$ ) with the normalized consumption rate. Thus, $\Psi$ measures not only the direct difference between $(f,S)$ and $(f^*,S^*)$, but also the deviation of the normalized age profile from its equilibrium shape, together with the deviation between normalized birth and consumption rates, and between the substrate level and the normalized consumption.

Estimate (3.25) of Theorem 1 guarantees global exponential convergence of the functional $\Psi(f,S)$ defined by (3.18) and Uniform Global Asymptotic Stability with respect to the measure $\Psi$. Definition (3.18) and inequality (2.9) show exponential convergence of $f[t]$ to $f^*$ in any $L^p$ norm with $p\in[1,+\infty]$. Since $d\left((f,S),(f^*,S^*)\right)\le\Psi(f,S)$ (recall (2.5) and (3.18)) estimate (3.25) shows global exponential convergence of the state to the equilibrium point in the metric of the state space. The global $KL$ stability estimate (3.25) and definition (3.18) of $\Psi$ can be exploited to conclude standard (uniform) global asymptotic stability properties with respect to various metrics.

**(b)** The sufficient conditions (3.20)-(3.24) that guarantee the global $KL$ stability estimate (3.25) are expressed in terms of the model data and auxiliary constants. The example in the following section illustrates their use and shows that the sufficient conditions (3.20)-(3.24) are not conservative.

We next provide a more detailed interpretation of these conditions:

First, condition (3.20) requires the quantity $2\langle q,r\rangle\alpha$ to dominate the auxiliary product $\omega\lambda$ where $\langle q,r\rangle$ is a weighted integral of the substrate consumption with respect to the normalized equilibrium profile $r(a)$ in (2.12). Since $\omega$ and $\lambda$ are auxiliary constants, their product can be chosen arbitrarily small and therefore, condition (3.20) can be satisfied without imposing any restrictive condition on the model parameters.

Next, conditions (3.21) and (3.23) constrain the admissible values of the weight $B$ appearing in the integral term of the Lyapunov functional measuring the deviation of the age distribution from the equilibrium. Condition (3.21) may impose a lower requirement on $B$, ensuring that deviations of the age distribution are penalized strongly enough in the Lyapunov estimate. The lower bound in (3.21) reflects the need for the age-distribution term to compensate for the additional coupling effects generated by the functions $h,p$. Thus, $B$ must be chosen sufficiently large so that deviations in the age profile are penalized strongly enough. On the other hand, condition (3.23) imposes an upper bound on $B$ preventing this weight from becoming too large. Thus, conditions (3.21) and (3.23) require the existence of an admissible range of values for $B$ ensuring a proper balance in the contribution of the age-distribution term to the Lyapunov functional.

More specifically, the parameter $\lambda \in [0,1]$ distributes certain coupling terms between these inequalities and the dissipation condition (3.22). When $\lambda = 0$, condition (3.21) becomes nonrestrictive, while the terms multiplied by $1-\lambda$ in (3.22) attain their maximal contribution. In this case, the functions $h(a)$ and $p(a)$, which describe additional terms in the differential relations for $q(a)$ and $k(a)$ in (3.2), contribute significantly to condition (3.22). Note that conditions (3.21) and (3.22) can be trivially satisfied when $h(a) \equiv 0$ and $p(a) \equiv 0$. Consequently, the parameters $\sigma$ and $\gamma$, together with the minimal mortality $L = \inf_{a \geq 0} \beta(a)$, must dominate the additional coupling effects generated by $h$ and $p$ in (3.22). For $\lambda > 0$, part of this influence is shifted to condition (3.21), thereby relaxing the dissipation requirement in (3.22). Overall, conditions (3.20)–(3.23) ensure that the parameters of the Lyapunov functional can be chosen so that the dissipative mechanisms dominate the coupling effects arising from reproduction and substrate consumption.

Finally, condition (3.24) requires a certain matrix $P$ to be positive definite. This matrix appears in the estimate of the Lyapunov derivative through a term involving a vector built from the deviations of the substrate concentration $S(t)$ and the auxiliary variables $\varphi(t)$ and $w(t)$. The positive definiteness of $P$ ensures that these deviations, taken together, are effectively controlled, so that no combination of $S$, $\varphi$, and $w$ can offset their contribution in the stability estimate.

## 4. Illustrative Example

In this section we apply Theorem 1 to the age-structured chemostat model studied in [32] and derive an explicit condition guaranteeing global asymptotic stability of the equilibrium. The model studied in [32] is given by

$$\frac{\partial f}{\partial t}(t,a) + \frac{\partial f}{\partial a}(t,a) = -(L+D) f(t,a) \tag{4.1}$$

$$f(t,0) = \mu(S) Y \int_0^{+\infty} \exp(-\tilde{k}a) f(t,a) da \tag{4.2}$$

$$\dot{S}(t) = D(S_{in} - S(t)) - \mu(S(t)) \int_0^{+\infty} f(t,a) da \tag{4.3}$$

where $Y>0$, $\tilde{k}\geq 0$, $L>0$, $D>0$ and $\mu\in C^1(\mathbb{R}_+)$, increasing with $\mu(0)=0$, and $\mu(S)>0$ for $S>0$. Model (4.1), (4.2), (4.3) corresponds to system (2.1), (2.2), (2.3) with

$$\beta(a)\equiv L,\quad q(a)\equiv 1,\quad k(a)=Y\exp(-\tilde{k}a) \tag{4.4}$$

In this case

$$r(a)=\exp(-(D+L)a) \tag{4.5}$$

while a straightforward calculation yields

$$\langle k,r\rangle=\frac{Y}{D+L+\tilde{k}}\qquad \langle q,r\rangle=\frac{1}{D+L} \tag{4.6}$$

The equilibrium is determined by the Lotka-Sharpe condition (2.10), which, due to (4.6), in this case is given by the point $S^*\in(0,S_{in})$ that satisfies

$$\mu\left(S^*\right)=\frac{D+L+\tilde{k}}{Y} \tag{4.7}$$

Using (4.6) and (2.11), the corresponding equilibrium for the age distribution is given by

$$f^*(a)=\frac{YD\left(D+L\right)\left(S_{in}-S^*\right)}{\left(D+L+\tilde{k}\right)}\exp\left(-(D+L)a\right) \tag{4.8}$$

**Proposition 1:** *Consider the age-structured chemostat model (4.1), (4.2), (4.3) with* $Y>0$, $\tilde{k}\geq 0$, $L>0$, $D>0$. *Let* $(f^*,S^*)$ *be the equilibrium defined by (4.7) and (4.8). If*

$$D>\frac{\tilde{k}^2}{8\left(2L+\tilde{k}\right)} \tag{4.9}$$

*then,* $(f^*,S^*)$ *is uniformly globally asymptotically stable with respect to the measure* $\Psi$ *defined by (3.18).*

**Remark on Proposition 1:** It was shown in [32], that after the conversion of system (4.1), (4.2), (4.3) to a system of ordinary differential equations, its linearization gives local attractivity of the equilibrium under the condition that

$$D\geq\frac{\left(\tilde{k}-L\right)^2}{8\tilde{k}} \tag{4.10}$$

On the other hand, condition (4.9) and estimate (3.25) establish (uniform) global asymptotic stability of the equilibrium (4.7), (4.8). Notice that when $3\tilde{k}>2L$, condition (4.9) is more conservative than (4.10), whereas, when $3\tilde{k}<2L$, condition (4.9) is less conservative than (4.10). Both conditions give the same limit condition $D>\tilde{k}/8$ as $L\to 0^+$ (when the mortality rate is negligible), which means that there must be roughly eight generations of organisms in the chemostat for instability of the equilibrium to be possible.

## 5. Proofs

In this section we provide the proofs of all results. We start with the proof of Lemma 1.

**Proof of Lemma 1:** Let arbitrary $(f_0, S_0) \in X$ and $D > 0$ be given and consider the weak solution $(f[t], S(t))$ with input $D$ of the initial-boundary value problem (2.1), (2.2), (2.3), (2.6). Let also $(f^*, S^*) \in X$ be an equilibrium point for (2.1), (2.2), (2.3) that corresponds to the constant input $D > 0$. Using the boundary condition (2.2) and that at equilibrium we have $f^*(a) = f^*(0)r(a)$, we get by Lemma 1 in [18] and (2.11), (2.12) that

$$f(t,a) - f^*(a) = \begin{cases} \left(f_0(a-t) - f^*(a-t)\right)\exp\left(-Dt - \int\limits_{a-t}^{a} \beta(s)ds\right) & for \quad 0 \le t \le a \\ \left(x(t-a) - f^*(0)\right)\exp\left(-Da - \int\limits_{0}^{a} \beta(s)ds\right) & for \quad t > a \ge 0 \end{cases} \tag{5.1}$$

where $x(t) = \mu\left(S(t)\right)\langle k, f[t]\rangle$.

Using (5.1) we get

$$\begin{aligned} \left\| f[t] - f^* \right\|_1 &= \int\limits_0^t \left| f(t,a) - f^*(a) \right| da + \int\limits_t^{+\infty} \left| f(t,a) - f^*(a) \right| da \\ &= \int\limits_0^t \left| x(t-a) - f^*(0) \right| \exp\left(-Da - \int\limits_0^a \beta(s)ds\right) da \\ &+ \exp(-Dt) \int\limits_t^{+\infty} \left| f_0(a-t) - f^*(a-t) \right| \exp\left(-\int\limits_{a-t}^{a} \beta(s)ds\right) da \\ &= \exp(-Dt) \int\limits_0^t \left| x(l) - f^*(0) \right| \exp\left(Dl - \int\limits_0^{t-l} \beta(s)ds\right) dl \\ &+ \exp(-Dt) \int\limits_0^{+\infty} \left| f_0(l) - f^*(l) \right| \exp\left(-\int\limits_l^{t+l} \beta(s)ds\right) dl \end{aligned} \tag{5.2}$$

Since $\beta(a) \ge 0$, $a \ge 0$, $D > 0$ we obtain from (5.2) the following estimate

$$\left\| f[t] - f^* \right\|_1 \le \int\limits_0^t \left| x(l) - f^*(0) \right| dl + \left\| f_0 - f^* \right\|_1 \tag{5.3}$$

Recall that $x(t)=\mu(S(t))\langle k,f[t]\rangle$ and $f^*(0)=\mu(S^*)\langle k,f^*\rangle$. Then, since $\mu$ is bounded and since $k\in L^\infty$ we get

$$\begin{aligned}&\left|x(t)-f^*(0)\right|=\left|\mu\left(S(t)\right)\langle k,f[t]\rangle-\mu\left(S^*\right)\left\langle k,f^*\right\rangle\right|\\&\le\mu\left(S(t)\right)\left|\left\langle k,f[t]-f^*\right\rangle\right|+\left\langle k,f^*\right\rangle\left|\mu\left(S(t)\right)-\mu\left(S^*\right)\right|\\&\le\mu_{\max}\left\|k\right\|_\infty\left\|f[t]-f^*\right\|_1+\left\langle k,f^*\right\rangle L_\mu\left|S(t)-S^*\right|\end{aligned}\tag{5.4}$$

where $L_\mu=\max\left\{\mu'(S):S\in\left[0,S_{in}\right]\right\}$ and $\mu_{\max}=\mu(S_{in})$ (recall that $\mu$ is increasing and $S(t)\in(0,S_{in})$ for all $t\ge 0$). Using (5.3) and (5.4) we obtain the following estimate

$$\left\|f[t]-f^*\right\|_1\le\mu_{\max}\left\|k\right\|_\infty\int_0^t\left\|f[l]-f^*\right\|_1dl+\left\|f_0-f^*\right\|_1+\left\langle k,f^*\right\rangle L_\mu\int_0^t\left|S(l)-S^*\right|dl\tag{5.5}$$

Since $(f^*,S^*)$ is an equilibrium of (2.1), (2.2), (2.3), equation (2.3) can also be written as

$$\dot S(t)=-D\left(S(t)-S^*\right)-\mu\left(S(t)\right)\langle q,f[t]\rangle+\mu\left(S^*\right)\left\langle q,f^*\right\rangle$$

which gives

$$\begin{aligned}&S(t)-S^*=\exp\left(-Dt\right)\left(S_0-S^*\right)\\&-\int_0^t\exp\left(-D(t-l)\right)\left(\mu\left(S(l)\right)\langle q,f[l]\rangle-\mu\left(S^*\right)\left\langle q,f^*\right\rangle\right)dl\end{aligned}\tag{5.6}$$

Using the fact that $q\in L^\infty$ and definitions $L_\mu=\max\left\{\mu'(S):S\in\left[0,S_{in}\right]\right\}$ and $\mu_{\max}=\mu(S_{in})$ we directly obtain from (5.6) that

$$\begin{aligned}&\left|S(t)-S^*\right|\le\left|S_0-S^*\right|+\int_0^t\left|\mu\left(S(l)\right)\langle q,f[l]\rangle-\mu\left(S^*\right)\left\langle q,f^*\right\rangle\right|dl\\&\le\left|S_0-S^*\right|+\int_0^t\mu\left(S(l)\right)\left|\left\langle q,f[l]-f^*\right\rangle\right|dl+\left\langle q,f^*\right\rangle\int_0^t\left|\mu\left(S(l)\right)-\mu\left(S^*\right)\right|dl\\&\le\left|S_0-S^*\right|+\mu_{\max}\left\|q\right\|_\infty\int_0^t\left\|f[l]-f^*\right\|_1dl+\left\langle q,f^*\right\rangle L_\mu\int_0^t\left|S(l)-S^*\right|dl\end{aligned}\tag{5.7}$$

Define $\bar R:=\max\left(\mu_{\max}\left(\left\|k\right\|_\infty+\left\|q\right\|_\infty\right),\left\langle q+k,f^*\right\rangle L_\mu\right)$. Estimates (5.5) and (5.7) give for $t\ge 0$ that

$$\left\| f[t]-f^{*} \right\|_{1} + \left| S(t)-S^{*} \right| \le \left\| f_{0}-f^{*} \right\|_{1} + \left| S_{0}-S^{*} \right| + \bar{R}\int_{0}^{t} \left( \left\| f[l]-f^{*} \right\|_{1} + \left| S(l)-S^{*} \right| \right) dl \tag{5.8}$$

Inequality (5.8) and a direct application of Gronwall-Bellman inequality imply that the following estimate holds for $t \ge 0$

$$\left\| f[t]-f^{*} \right\|_{1} + \left| S(t)-S^{*} \right| \le \exp\left(\bar{R}\,t\right)\left( \left\| f_{0}-f^{*} \right\|_{1} + \left| S_{0}-S^{*} \right| \right) \tag{5.9}$$

Inequality (2.17) follows immediately from definition (2.5) and (5.9).

Next, from formula (5.1), (5.4), (5.9) and the fact that $\beta(a) \ge 0$, $a \ge 0$, and $D > 0$ we get

$$\begin{aligned}
&\left\| f[t]-f^{*} \right\|_{\infty} \le \max\left( \left\| f_{0}-f^{*} \right\|_{\infty}, \max_{l\in[0,t]}\left( \left| x(l)-f^{*}(0) \right| \right) \right) \\
&\le \max\left( \left\| f_{0}-f^{*} \right\|_{\infty}, \max\left( \mu_{\max}\left\| k \right\|_{\infty}, \left\langle k, f^{*} \right\rangle L_{\mu} \right) \max_{l\in[0,t]}\left( \left\| f[l]-f^{*} \right\|_{1} + \left| S(l)-S^{*} \right| \right) \right) \\
&\le \max\left( \left\| f_{0}-f^{*} \right\|_{\infty}, \max\left( \mu_{\max}\left\| k \right\|_{\infty}, \left\langle k, f^{*} \right\rangle L_{\mu} \right) \exp\left(\bar{R}\,t\right)\left( \left\| f_{0}-f^{*} \right\|_{1} + \left| S_{0}-S^{*} \right| \right) \right)
\end{aligned} \tag{5.10}$$

Inequality (2.18) is a consequence of (5.10) and definition (2.5) with $\bar{C} = \max\left( \mu_{\max}\left\| k \right\|_{\infty}, \left\langle k, f^{*} \right\rangle L_{\mu} \right)$. Since $\mu(S) \ge 0$ for $S \ge 0$, and $\left\langle k, f \right\rangle > 0$ for all $\left(f, S\right) \in X$, equation (2.2) gives $\dot{S}(t) \le D\left(S_{in} - S(t)\right)$, $t \ge 0$, from which (2.19) follows directly. The proof is complete. ⊲

We next provide the proof of Lemma 2.

**Proof of Lemma 2:** Let $\left(f_{0}, S_{0}\right) \in X$. By virtue of Theorem 1 in [18], the weak solution $\left(f[t], S(t)\right) \in X$ of the initial-boundary value problem (2.1), (2.2), (2.3), (2.6) is unique and is defined for all $t \ge 0$.

Define

$$Y(t) := RS(t) + \left\| f[t] \right\|_{1}, \text{ for } t \ge 0 \tag{5.11}$$

Using (2.3), (3.1), (2.16), definition (5.11) and the fact that $f(t,a) > 0$ and $\beta(a) \ge 0$ we get for all $t \ge 0$

$$\begin{aligned}
&\dot{Y}(t) = D\left(RS_{in} - RS(t)\right) - \mu\left(S(t)\right)\left\langle Rq, f \right\rangle - \left\langle \beta, f \right\rangle - D\left\| f \right\|_{1} + \mu\left(S\right)\left\langle k, f \right\rangle \\
&\le D\left(RS_{in} - Y(t)\right) + \mu\left(S(t)\right)\left\langle k - Rq, f \right\rangle \\
&\le D\left(RS_{in} - Y(t)\right)
\end{aligned} \tag{5.12}$$

The differential inequality (5.12) gives the following estimate for all $t \ge 0$:

$$Y(t) \le RS_{in} + \left(Y(0) - RS_{in}\right)\exp\left(-Dt\right) \tag{5.13}$$

Define:

$$T_1\left(f_0, S_0\right) := \frac{1}{D}\ln\left(1 + \frac{\left(RS_0 + \|f_0\|_1 - F\right)^+}{F - RS_{in}}\right) \tag{5.14}$$

Then (5.13) and (5.14) guarantee that $Y(t) = RS(t) + \|f[t]\|_1 \le F$ for all $t \ge T_1\left(f_0, S_0\right)$. Moreover, if $Y(0) = RS_0 + \|f_0\|_1 \le F$ then $Y(t) = RS(t) + \|f[t]\|_1 \le F$ for all $t \ge 0$.

Since $\|f[t]\|_1 \le Y(t)$ (recall (5.11)), it follows that

$$\|f[t]\|_1 \le F \text{ for all } t \ge T_1\left(f_0, S_0\right) \tag{5.15}$$

Using the facts that $S(t) \in (0, S_{in})$ for all $t \ge 0$, (5.15), $\sup\left\{\frac{\mu(S)}{S} : S \in \left(0, S_{in}\right)\right\} \le L_\mu$ (a consequence of (3.5) and the fact that $\mu\left(0\right) = 0$), (3.4) and boundedness of $q$, we get for every $t \ge T_1\left(f_0, S_0\right)$:

$$\begin{aligned} \dot{S}(t) &= D\left(S_{in} - S(t)\right) - \mu\left(S(t)\right)\langle q, f\rangle \\ &\ge DS_{in} - DS(t) - \|q\|_\infty \mu\left(S(t)\right)\|f[t]\|_1 \\ &\ge DS_{in} - DS(t) - F\|q\|_\infty \mu\left(S(t)\right) \\ &\ge DS_{in} - \left(D + L_\mu F\|q\|_\infty\right)S(t) \\ &= \left(D + L_\mu F\|q\|_\infty\right)\left(2\underline{S} - S(t)\right) \end{aligned} \tag{5.16}$$

Then, from (5.16) we obtain the following estimate for $t \ge T_1\left(f_0, S_0\right)$:

$$S(t) \ge 2\underline{S} + \left(S\left(T_1\left(f_0, S_0\right)\right) - 2\underline{S}\right)\exp\left(-\left(D + L_\mu F\|q\|_\infty\right)\left(t - T_1\left(f_0, S_0\right)\right)\right) \tag{5.17}$$

Estimate (5.17) guarantees that if $RS_0 + \|f_0\|_1 \le F$ (which implies $T_1\left(f_0, S_0\right) = 0$; recall (5.14)) and $S(0) = S_0 \ge \underline{S}$ then $Y(t) = RS(t) + \|f[t]\|_1 \le F$ and $S(t) \ge \underline{S}$ for all $t \ge 0$. Estimate (5.17) guarantees that $S(t) \ge \underline{S}$ for all $t \ge T_2\left(f_0\right)$, where

$$T_2\left(f_0, S_0\right) := T_1\left(f_0, S_0\right) + \frac{\ln(2)}{D + L_\mu F\|q\|_\infty} \tag{5.18}$$

Combining (5.18) and (5.14), we conclude that (3.3) holds with $T\left(s\right) := \frac{1}{D}\ln\left(1 + \frac{\left(s - F\right)^+}{F - RS_{in}}\right) + \frac{\ln(2)}{D + L_\mu F\|q\|_\infty}$ for $s \ge 0$. The proof is complete. ◁

The proof of Theorem 1 requires certain auxiliary results. The following lemma provides additional estimates for the solutions of (2.1), (2.2), (2.3).

**Lemma 3:** *Suppose that assumptions (A), (B) hold. Then* $b \le 1$, $\gamma \le 1$ *and for every* $(f_0, S_0) \in X$ *the weak solution* $(f[t], S(t))$ *of the initial-boundary value problem (2.1), (2.2), (2.3), (2.6) satisfies the following estimates for all* $t \ge 0$

$$\xi(t) \ge \exp\left(D(b-1)t\right)\xi(0) \tag{5.19}$$

$$\zeta(t) \ge \exp\left(D(\gamma-1)t\right)\zeta(0) \tag{5.20}$$

*and*

$$S(t) \ge \min\left( S_0, \frac{DS_{in}}{D + L_\mu \|q\|_\infty \left( RS_{in} + \left( RS_0 + \|f_0\|_1 - RS_{in} \right)^+ \right)} \right) \tag{5.21}$$

**Proof of Lemma 3:** Equations (2.14), (2.15), (3.2) and definitions (3.7) give for every weak solution $(f[t], S(t))$ of (2.1), (2.2), (2.3):

$$\begin{aligned} \dot{\xi} &= \left(k(0)\mu(S) + bD - D\right)\xi + \frac{\alpha\langle q, f^* \rangle}{\theta\langle k, f^* \rangle}\zeta + \frac{\langle p, f \rangle}{\langle k, f^* \rangle} \\ \dot{\zeta} &= \left(q(0)\mu(S) + \delta\theta\right)\frac{\langle k, f^* \rangle}{\langle q, f^* \rangle}\xi + D(\gamma-1)\zeta + \frac{\langle h, f \rangle}{\langle q, f^* \rangle} \end{aligned} \tag{5.22}$$

At equilibrium, we obtain from (5.22) that

$$\begin{aligned} D(b-1) &= -\frac{\alpha\langle q, f^* \rangle}{\theta\langle k, f^* \rangle} - \frac{\langle p, f^* \rangle}{\langle k, f^* \rangle} - k(0)\mu\left(S^*\right) \\ D(\gamma-1) &= -\frac{\langle h, f^* \rangle}{\langle q, f^* \rangle} - \left(q(0)\mu\left(S^*\right) + \delta\theta\right)\frac{\langle k, f^* \rangle}{\langle q, f^* \rangle} \end{aligned}$$

The above equations together with the facts that $k, q, h, p, \mu : \mathbb{R}_+ \to \mathbb{R}_+$ are non-negative functions and $\theta > 0$, $\alpha, \delta \ge 0$ guarantee that $b \le 1$ and $\gamma \le 1$.

Moreover, equations (5.22) guarantee the following differential inequalities:

$$\begin{aligned} \dot{\xi} &\ge D(b-1)\xi \\ \dot{\zeta} &\ge D(\gamma-1)\zeta \end{aligned} \tag{5.23}$$

Estimates (5.19) and (5.20) are consequences of the differential inequalities (5.23).

Equation (5.13) shows that the following estimate holds for all $t \geq 0$:

$$\|f[t]\|_1 \leq RS_{in} + \left(RS_0 + \|f_0\|_1 - RS_{in}\right)^+ \tag{5.24}$$

Using (2.3), (3.5), and (5.24) we get for all $t \geq 0$:

$$\begin{aligned} \dot{S}(t) &= D\left(S_{in} - S(t)\right) - \mu\left(S(t)\right)\langle q, f\rangle \\ &\geq DS_{in} - DS(t) - \|q\|_\infty \mu\left(S(t)\right)\|f[t]\|_1 \\ &\geq DS_{in} - DS(t) - L_\mu \|q\|_\infty \|f[t]\|_1 S(t) \\ &\geq DS_{in} - \left(D + L_\mu \|q\|_\infty \left(RS_{in} + \left(RS_0 + \|f_0\|_1 - RS_{in}\right)^+\right)\right)S(t) \end{aligned} \tag{5.25}$$

Differential inequality (5.25) implies the following estimate for all $t \geq 0$:

$$S(t) \geq S_0 \exp\left(-\left(D + L_\mu \|q\|_\infty \left(RS_{in} + \left(RS_0 + \|f_0\|_1 - RS_{in}\right)^+\right)\right)t\right)$$
$$+\frac{DS_{in}\left(1 - \exp\left(-\left(D + L_\mu \|q\|_\infty \left(RS_{in} + \left(RS_0 + \|f_0\|_1 - RS_{in}\right)^+\right)\right)t\right)\right)}{D + L_\mu \|q\|_\infty \left(RS_{in} + \left(RS_0 + \|f_0\|_1 - RS_{in}\right)^+\right)}$$

Estimate (5.21) is obtained from the above estimate. The proof is complete. $\lhd$

Next, we give some useful lemmas for the functional $V$ defined by (3.15). Lemma 4 below shows first that $V$ defined by (3.15) is differentiable along the weak solutions of (2.1), (2.2), (2.3), (2.6) with constant input $D$ and provides the exact derivative formula of the Lyapunov functional.

**Lemma 4:** *Suppose that assumption (B) holds. Let arbitrary $\left(f_0, S_0\right) \in X$ be given and consider the weak solution $\left(f[t], S(t)\right)$ with input $D$ of the initial-boundary value problem (2.1), (2.2), (2.3), (2.6). Then the map $t \mapsto V(f[t], S(t))$, $t \geq 0$ is differentiable and the derivative of the Lyapunov functional $V(f, S)$ along the weak solution $\left(f[t], S(t)\right)$, $t \geq 0$, satisfies the following equation for $t \geq 0$:*

$$\begin{aligned} &\frac{d}{dt}\left(V(f[t], S(t))\right) = \left(\exp\left(\varphi(t)\right) - 1\right)\dot{\varphi}(t) + \Gamma\left(\exp\left(w(t)\right) - 1\right)\dot{w}(t) \\ &-B\int_0^{+\infty} \rho^2(a)\chi(t,a)\left(\left(\beta(a) + D + \frac{\dot{\zeta}(t)}{\zeta(t)}\right)\left(\chi(t,a) + r(a)\right) + r'(a)\right)da \\ &+Q'(S(t))\dot{S}(t) + \frac{B}{2}\chi^2(t,0) - \sigma B\|\rho\chi\|_2^2 \end{aligned} \tag{5.26}$$

**Proof of Lemma 4:** Let arbitrary $(f_0,S_0)\in X$ be given and consider the weak solution $(f[t],S(t))$ with constant input $D$ of the initial-boundary value problem (2.1), (2.2), (2.3), (2.6). Define $x(t)=\mu(S(t))\langle k,f[t]\rangle$, $t\geq 0$. Using the boundary condition (2.2) and that at equilibrium we have $f^*(a)=f^*(0)r(a)$, we get by Lemma 1 in [18] and (2.11), (2.12) that

$$f(t,a)-f^*(a)=\begin{cases}\left(f_0(a-t)-f^*(a-t)\right)\dfrac{r(a)}{r(a-t)} & \text{for} \quad 0\leq t\leq a\\ \left(x(t-a)-f^*(0)\right)r(a) & \text{for} \quad t>a\geq 0\end{cases} \tag{5.27}$$

Next, recalling from (3.8) that $\chi(a)=\dfrac{f(a)}{f^*(0)\zeta}-r(a)$ the Lyapunov functional in (3.15) can be written as

$$\begin{aligned}V(f,S)&=\exp(\varphi)-\varphi-1+\Gamma\left(\exp(w)-w-1\right)+Q(S)\\&+\frac{B}{2\zeta^2\left(f^*(0)\right)^2}\int_0^{+\infty}\rho^2(a)\left(f(a)-f^*(a)\right)^2da\\&+\frac{B(\zeta-1)^2}{2\zeta^2}\int_0^{+\infty}\rho^2(a)r^2(a)da-\frac{B(\zeta-1)}{\zeta^2 f^*(0)}\int_0^{+\infty}\rho^2(a)r(a)\left(f(a)-f^*(a)\right)da\end{aligned} \tag{5.28}$$

Splitting the first and last integral of (5.28) at $a=t$, we get, via (5.27) and a change of variables $\tau=t-a$, that

$$\begin{aligned}\int_0^{+\infty}\rho^2(a)r(a)\left(f(t,a)-f^*(a)\right)da&=\int_0^{t}\rho^2(t-s)r^2(t-s)\left(x(s)-f^*(0)\right)ds\\&+\int_0^{+\infty}\rho^2(t+\tau)\frac{r^2(t+\tau)}{r(\tau)}\left(f_0(\tau)-f^*(\tau)\right)d\tau\end{aligned} \tag{5.29}$$

$$\begin{aligned}\int_0^{+\infty}\rho^2(a)\left(f(a)-f^*(a)\right)^2da&=\int_0^{t}\rho^2(t-s)r^2(t-s)\left(x(s)-f^*(0)\right)^2ds\\&+\int_0^{+\infty}\rho^2(t+\tau)\frac{r^2(t+\tau)}{r^2(\tau)}\left(f_0(\tau)-f^*(\tau)\right)^2d\tau\end{aligned} \tag{5.30}$$

Combining (5.28), (5.29), and (5.30) we get

$$V(f[t],S(t))=\exp\left(\varphi(t)\right)-\varphi(t)-1+\Gamma\left(\exp\left(w(t)\right)-w(t)-1\right)+Q\left(S(t)\right)$$

$$+\frac{B}{2\zeta^2(t)\left(f^*(0)\right)^2}\int_0^t \rho^2(t-s)r^2(t-s)\left(x(s)-f^*(0)\right)^2 ds$$

$$+\frac{B}{2\zeta^2(t)\left(f^*(0)\right)^2}\int_0^{+\infty} \rho^2(t+\tau)\frac{r^2(t+\tau)}{r^2(\tau)}\left(f_0(\tau)-f^*(\tau)\right)^2 d\tau$$

$$+\frac{B\left(\zeta(t)-1\right)^2}{2\zeta^2(t)}\int_0^{+\infty} \rho^2(a)r^2(a)da-\frac{B\left(\zeta(t)-1\right)}{\zeta^2(t)f^*(0)}\int_0^t \rho^2(t-s)r^2(t-s)\left(x(s)-f^*(0)\right)ds$$

$$-\frac{B\left(\zeta(t)-1\right)}{\zeta^2(t)f^*(0)}\int_0^{+\infty} \rho^2(t+\tau)\frac{r^2(t+\tau)}{r(\tau)}\left(f_0(\tau)-f^*(\tau)\right)d\tau$$

Using Leibniz's integral rule and equations (2.11), (2.12), (3.8), (3.11), (3.16), (5.29) and (5.30) we obtain (5.26). The proof is complete. $\triangleleft$

Next, we give a technical lemma whose proof can be found in the appendix. We use the notation $\nabla^2 W(0)$ for the Hessian matrix at $0\in O\subseteq\mathbb{R}^n$ of a function of $W\in C^2(O)$, where $O$ is an open set.

**Lemma 5:** *Let $O\subseteq\mathbb{R}^n$ be an open set with $0\in O$. Suppose that $W\in C^2\left(O\right)$ is a positive definite function with $\nabla^2 W(0)$ being positive definite and such that the set $\left\{x\in O: W(x)\le r\right\}$ is compact for all $r\ge 0$. Then there exists a non-decreasing function $\bar{p}:\mathbb{R}_+\to\mathbb{R}_+$ such that $W(x)\bar{p}\left(W(x)\right)\ge |x|^2$ for all $x\in O$.*

The following lemma exploits Lemma 5 and establishes some useful estimates for the Lyapunov functional $V$ in (3.15).

**Lemma 6:** *There exist non-decreasing functions $\bar{\rho}:\mathbb{R}_+\to\mathbb{R}_+$, $\bar{\gamma}:\mathbb{R}_+\to(0,+\infty)$ such that the following inequalities holds for all $\left(f,S\right)\in X$ :*

$$\left|\mu\left(S\right)\left\langle k,f\right\rangle-f^*(0)\right|+\left|S-S^*\right|\le\left(V\left(f,S\right)\right)^{1/2}\bar{\rho}\left(V\left(f,S\right)\right) \tag{5.31}$$

$$V\left(f,S\right)\le\bar{\gamma}\left(V\left(f,S\right)\right)E(f,S) \tag{5.32}$$

*where*

$$E(f,S)=\left(\exp\left(\varphi\right)-1\right)^2+\left(\exp\left(w\right)-1\right)^2+\left(g(S)-1\right)^2+\left\|\rho\chi\right\|_2^2 \tag{5.33}$$

**Proof of Lemma 6:** Let $\varphi, w$ as defined by (3.9). Notice that from (2.2) we get the equilibrium relation $f^*(0)=\mu(S^*)\left\langle k,f^*\right\rangle$. Using definitions (3.9), (3.10), and the previous equality, we get

$$\mu(S)\langle k,f\rangle - f^*(0) = f^*(0)\left(\frac{(S_{in}-S)g(S)}{(S_{in}-S^*)}\exp(\varphi+w)-1\right) \tag{5.34}$$

By adding and subtracting terms, (5.34) gives

$$\begin{aligned}
\mu(S)\langle k,f\rangle - f^*(0) &= f^*(0)\frac{(S_{in}-S)g(S)}{(S_{in}-S^*)}\left(\exp(\varphi)-1\right) \\
&+ f^*(0)\frac{(S_{in}-S)g(S)}{(S_{in}-S^*)}\left(\exp(w)-1\right)\left(\exp(\varphi)-1\right) \\
&+ f^*(0)\frac{(S_{in}-S)g(S)}{(S_{in}-S^*)}\left(\exp(w)-1\right) \\
&+ \frac{f^*(0)(S_{in}-S)}{(S_{in}-S^*)\mu(S^*)}\left(\mu(S)-\mu(S^*)\right) + \frac{f^*(0)}{(S_{in}-S^*)}\left(S^*-S\right)
\end{aligned} \tag{5.35}$$

For $S \in (0, S_{in})$, by using the fact that $g(S)$ is increasing (since $\mu$ is increasing) and definition (3.5), we get

$$\begin{aligned}
\left|\mu(S)\langle k,f\rangle - f^*(0)\right| &\le f^*(0)\frac{S_{in}g(S_{in})}{(S_{in}-S^*)}\left|\exp(\varphi)-1\right| \\
&+ f^*(0)\frac{S_{in}g(S_{in})}{(S_{in}-S^*)}\left|\exp(w)-1\right|\left|\exp(\varphi)-1\right| \\
&+ f^*(0)\frac{S_{in}g(S_{in})}{(S_{in}-S^*)}\left|\exp(w)-1\right| + \frac{f^*(0)}{(S_{in}-S^*)}\left(\frac{S_{in}L_\mu}{\mu(S^*)}+1\right)\left|S-S^*\right|
\end{aligned} \tag{5.36}$$

Due to Lemma 5, there exists a non-decreasing function $\tilde{p}:\mathbb{R}_+ \to \mathbb{R}_+$ such that $|y|^2 \le \left(\exp(y)-y-1\right)\tilde{p}\left(\exp(y)-y-1\right)$ for all $y \in \mathbb{R}$. Using the inequality $\left|\exp(y)-1\right| \le |y|\exp\left(|y|\right)$ that holds for all $y \in \mathbb{R}$, we get

$$\left|\exp(y)-1\right| \le \left(\exp(y)-y-1\right)^{1/2}\varpi\left(\exp(y)-y-1\right) \tag{5.37}$$

for all $y \in \mathbb{R}$ with $\varpi(s) := \left(\tilde{p}(s)\right)^{1/2}\exp\left(s^{1/2}\left(\tilde{p}(s)\right)^{1/2}\right)$. Notice that $\varpi:\mathbb{R}_+ \to \mathbb{R}_+$ is a non-decreasing function.

Using now (3.15) and (5.37) with $y=\varphi$ and $y=w$ we get, respectively,

$$\left|\exp(\varphi)-1\right| \le \left(V(f,S)\right)^{1/2}\varpi\left(V(f,S)\right) \tag{5.38}$$

$$\left|\exp(w)-1\right| \le \Gamma^{-1/2}\left(V(f,S)\right)^{1/2} \varpi\left(\Gamma^{-1}V(f,S)\right) \tag{5.39}$$

Define for $x \in \left(-S^*, S_{in}-S^*\right)$:

$$\ell(x) := Q\left(S^*+x\right) \tag{5.40}$$

where $Q(S)$ satisfies (3.17). Then, using (3.17) and (5.40) we get that $\ell(0)=0$ and

$$\ell(x) = M \int_{S^*}^{S^*+x} \frac{\left(g(l)-1\right)}{(S_{in}-l)\,g(l)} dl \tag{5.41}$$

Notice now that

$$\ell'(x) = M \frac{\left(g\left(S^*+x\right)-1\right)}{\left(S_{in}-S^*-x\right) g\left(S^*+x\right)}$$

$$\ell''(x) = M \frac{g'\left(S^*+x\right)\left(S_{in}-S^*-x\right)+\left(g\left(S^*+x\right)-1\right) g\left(S^*+x\right)}{\left(S_{in}-S^*-x\right)^2 g^2\left(S^*+x\right)}$$

from which, the following properties are obtained

$$\begin{gathered} x \ne 0 \Rightarrow x\ell'(x) > 0 \\ \ell''(0) = M \frac{g'\left(S^*\right)}{\left(S_{in}-S^*\right)} > 0 \end{gathered} \tag{5.42}$$

where we have used the fact that $g(S)$ is increasing with $g\left(S^*\right)=1$.

Next, we show that the set

$$\left\{x \in \left(-S^*, S_{in}-S^*\right) : \ell(x) \le r\right\} \tag{5.43}$$

is compact for all $r \ge 0$. Since $g$ is increasing, we can define $\tilde{L} = \max\left\{ g'(S) : S \in \left[0, S^*\right]\right\}$, which gives for $S \in \left[0, S^*\right]$ that $g(S) \le \tilde{L}S$.

For $-S^* < x \le -\frac{S^*}{2}$ and by using the fact that $g$ is increasing, inequalities $g(S) \le \tilde{L}S$, $1-g(S) \ge 0$, $S \in \left[0, S^*\right]$ and (5.41), we get

$$\begin{aligned}\ell(x) &\geq \frac{M}{\tilde{L}} \int_{S^*+x}^{S^*} \frac{(1-g(l))}{(S_{in}-l)l} dl \geq \frac{M}{\tilde{L}S_{in}} \int_{S^*+x}^{S^*} \frac{(1-g(l))}{l} dl \\ &\geq \frac{M}{\tilde{L}S_{in}} \int_{S^*+x}^{S^*/2} \frac{(1-g(l))}{l} dl \geq \frac{M\left(1-g\left(S^*/2\right)\right)}{\tilde{L}S_{in}} \int_{S^*+x}^{S^*/2} \frac{1}{l} dl \\ &\geq \frac{M\left(1-g\left(S^*/2\right)\right)}{\tilde{L}S_{in}} \ln\left(\frac{S^*}{2\left(S^*+x\right)}\right)\end{aligned} \tag{5.44}$$

For $S_{in} - S^* > x \geq \frac{S_{in} - S^*}{2}$ and by using similar arguments as above we also get that

$$\begin{aligned}\ell(x) &\geq \frac{M}{g(S_{in})} \int_{S^*}^{S^*+x} \frac{(g(l)-1)}{(S_{in}-l)} dl \geq \frac{M}{g(S_{in})} \int_{(S_{in}+S^*)/2}^{S^*+x} \frac{(g(l)-1)}{(S_{in}-l)} dl \\ &\geq \frac{M\left(g\left(\left(S_{in}+S^*\right)/2\right)-1\right)}{g(S_{in})} \int_{(S_{in}+S^*)/2}^{S^*+x} \frac{dl}{(S_{in}-l)} \\ &\geq \frac{M\left(g\left(\left(S_{in}+S^*\right)/2\right)-1\right)}{g(S_{in})} \ln\left(\frac{\left(S_{in}-S^*\right)}{2\left(S_{in}-S^*-x\right)}\right)\end{aligned} \tag{5.45}$$

Inequalities (5.44) and (5.45) show that $\lim_{x\to(-S^*)^+}\left(\ell(x)\right) = +\infty$ and $\lim_{x\to(S_{in}-S^*)^-}\left(\ell(x)\right) = +\infty$, which imply that the set in (5.43) is compact for all $r \geq 0$.

Since the set in (5.43) is compact for all $r \geq 0$ and due to (5.42), it follows by Lemma 5, that there exists a non-decreasing function $\bar{p}: \mathbb{R}_+ \to \mathbb{R}_+$ such that $|x|^2 \leq \ell(x)\,\bar{p}(\ell(x))$ for all $x \in \left(-S^*, S_{in} - S^*\right)$. Thus, for $x = S - S^*$, definition (5.40) gives

$$\left|S - S^*\right|^2 \leq Q(S)\,\bar{p}(Q(S)) \text{ for all } S \in \left(0, S_{in}\right) \tag{5.46}$$

Now define $\tilde{\gamma}(s) := \left(\bar{p}(s)\right)^{1/2}$, $s \geq 0$. Since $\bar{p}$ is nondecreasing, so is $\tilde{\gamma}$. Thus, since by (3.15) we have that $Q(S) \leq V(f,S)$, $(f,S) \in X$, it follows that

$$\left|S - S^*\right| \leq \left(V(f,S)\right)^{1/2} \tilde{\gamma}\left(V(f,S)\right) \text{ for all } S \in \left(0, S_{in}\right) \tag{5.47}$$

Combining (5.36), (5.38), (5.39), and (5.47), we finally obtain (5.31) with

$$\bar{\rho}(r) := f^*(0)\frac{S_{in}g\left(S_{in}\right)}{\left(S_{in}-S^*\right)}\varpi\left(r\right)\left(1+\frac{r^{1/2}}{\Gamma^{1/2}}\varpi\left(\Gamma^{-1}r\right)\right)+\tilde{\gamma}\left(r\right)$$

$$+\frac{f^*(0)}{\left(S_{in}-S^*\right)}\left(\frac{S_{in}g\left(S_{in}\right)}{\Gamma^{1/2}}\varpi\left(\Gamma^{-1}r\right)+\left(\frac{S_{in}L_{\mu}}{\mu\left(S^*\right)}+1\right)\tilde{\gamma}\left(r\right)\right)$$

Next, we show that (5.32) holds for some non-decreasing function $\bar{\gamma}$. Define

$$\Theta(S) = \begin{cases} \dfrac{Q(S)}{(g(S)-1)^2} & S \in \left(0,S^*\right) \cup \left(S^*,S_{in}\right) \\ \dfrac{M}{2\left(S_{in}-S^*\right)g'(S^*)} & S = S^* \end{cases} \tag{5.48}$$

which due to (3.17) is continuous since $g,Q \in C^1$ with $g(S^*)=1$ and $Q(S^*)=0$. Define also for $s \geq 0$

$$\gamma_Q(s) := \sup\left\{\Theta(S) : S \in (0,S_{in}), Q(S) \leq s\right\}$$

which is non-decreasing and well defined for each $s \geq 0$, by virtue of (5.40) and (5.43). Since $Q(S) \leq V(f,S)$ for all $(f,S) \in X$, it follows that

$$Q(S) \leq \gamma_Q(V(f,S))\left(g(S)-1\right)^2 \tag{5.49}$$

Define for $s > 0$

$$\beta_1(s) := \sup\left(\frac{\exp(y)-y-1}{\left(\exp(y)-1\right)^2} : y \neq 0, e^y - y - 1 \leq s\right) \tag{5.50}$$

which is well defined (due to the fact that the set $\left\{y \in \mathbb{R} : \exp(y)-y-1 \leq s\right\}$ is compact for all $s \geq 0$) and non-decreasing. For any $y \in \mathbb{R}$ with $\exp(y)-y-1 \leq s$, $s \geq 0$, definition (5.50) implies that $\exp(y)-y-1 \leq \beta_1(s)\left(\exp(y)-1\right)^2$. Thus,

$$\exp\left(\varphi\right)-\varphi-1 \leq \beta_1(V(f,S))\left(\exp(\varphi)-1\right)^2 \tag{5.51}$$

$$\exp\left(w\right)-w-1 \leq \beta_1(V(f,S)/\Gamma)\left(\exp(w)-1\right)^2 \tag{5.52}$$

Then due to (5.51), (5.52), (5.49) and (3.15) we obtain (5.32) with

$$\bar{\gamma}(s) := \max\left(\beta_1(s), \Gamma\beta_1(s/\Gamma), \gamma_Q(s), B/2\right)$$

The proof is complete. $\triangleleft$

The following lemma establishes that the exponential decay of the Lyapunov functional inside the trapping region $\Omega$ defined by (3.6) implies exponential convergence of the state to the equilibrium $(f^*, S^*)$ in terms of the metric $d$ defined by (2.5) and in terms of the sup-norm of $f$.

**Lemma 7:** *Suppose that assumptions (A), (B) hold. Suppose that there exist a non-increasing function $\bar{\sigma}:\mathbb{R}_+ \to (0,+\infty)$ and a functional $H : X \to \mathbb{R}_+$ such that for every $(f_0, S_0) \in X$ the weak solution of the initial-boundary value problem (2.1), (2.2), (2.3), (2.6) satisfies the following estimate for all $t \geq 0$ with $\bar{\sigma} := \bar{\sigma}\left(H\left(f_0, S_0\right)\right)$:*

$$V\left(f[t], S(t)\right) \leq \exp\left(-4\bar{\sigma}\left(H\left(f_0, S_0\right)\right)t\right) H\left(f_0, S_0\right) \tag{5.53}$$

*Then there exists a function $\kappa \in K_\infty$ such that for every $(f_0, S_0) \in X$ the weak solution of the initial-boundary value problem (2.1), (2.2), (2.3), (2.6) satisfies the following estimates for all $t \geq 0$:*

$$\left|S(t)-S^*\right| + \left\|f[t]-f^*\right\|_1 \leq \exp\left(-\bar{\sigma}\left(H\left(f_0, S_0\right)\right)t\right)\left(\left\|f_0 - f^*\right\|_1 + \kappa\left(H\left(f_0, S_0\right)\right)\right) \tag{5.54}$$

$$\left\|f[t]-f^*\right\|_\infty \leq \exp\left(-\bar{\sigma}\left(H\left(f_0, S_0\right)\right)t\right)\max\left(\left\|f_0 - f^*\right\|_\infty, \kappa\left(H\left(f_0, S_0\right)\right)\right) \tag{5.55}$$

**Proof of Lemma 7:** Without loss of generality we next assume that $\bar{\sigma}(s) \leq D/2$ for all $s \geq 0$.

Let $(f_0, S_0) \in X$ be given and consider the weak solution $\left(f[t], S(t)\right)$ of the initial-boundary value problem (2.1), (2.2), (2.3), (2.6). In what follows we use the notation $\bar{\sigma} := \bar{\sigma}\left(H\left(f_0, S_0\right)\right)$.

Define $x(t) = \mu\left(S(t)\right)\langle k, f[t]\rangle$, $t \geq 0$. Using the boundary condition (2.2) and that at equilibrium we have $f^*(a) = f^*(0) r(a)$, we get by Lemma 1 in [18] and (2.11), (2.12) that

$$f(t,a) - f^*(a) = \begin{cases} \left(f_0(a-t) - f^*(a-t)\right)\exp\left(-Dt - \int\limits_{a-t}^{a} \beta(s)ds\right) & \text{for} \quad 0 \leq t \leq a \\ \left(x(t-a) - f^*(0)\right)\exp\left(-Da - \int\limits_{0}^{a} \beta(s)ds\right) & \text{for} \quad t > a \geq 0 \end{cases} \tag{5.56}$$

Notice now that since $\beta(a) \geq 0$ for $a \geq 0$ and $D > 0$, we have that $\exp\left(-Da - \int_0^a \beta(s)ds\right) \leq 1$ and $\exp\left(-\int_{a-t}^{a} \left(\beta(s) + D\right)ds\right) \leq 1$. Consequently, from (5.56) we get that the following estimates hold for $t \geq 0$:

$$\begin{aligned}\left\|f[t]-f^*\right\|_1 &= \int_0^t \left|f(t,a)-f^*(a)\right|da + \int_t^{+\infty} \left|f(t,a)-f^*(a)\right|da \\ &\leq \int_0^t \left|x(l)-f^*(0)\right|dl + \exp\left(-Dt\right)\left\|f_0-f^*\right\|_1\end{aligned} \tag{5.57}$$

and

$$\left|f(t,a)-f^*(a)\right| \leq \begin{cases} \exp\left(-Dt\right)\left\|f_0-f^*\right\|_\infty & 0 \leq t \leq a \\ \left|x(t-a)-f^*(0)\right| & t > a \geq 0 \end{cases} \tag{5.58}$$

Due to continuity of $x(t)$, estimate (5.58) implies that for $t \geq 0$:

$$\left\|f[t]-f^*\right\|_\infty \leq \max\left(\exp\left(-Dt\right)\left\|f_0-f^*\right\|_\infty, \max_{l\in[0,t]}\left(\left|x(l)-f^*(0)\right|\right)\right) \tag{5.59}$$

Let arbitrary $t \geq 0$ be given. Estimates (5.57) and (5.59) give at time $t/2$:

$$\left\|f[t/2]-f^*\right\|_1 \leq \int_0^{t/2} \left|x(l)-f^*(0)\right|dl + \exp\left(-Dt/2\right)\left\|f_0-f^*\right\|_1 \tag{5.60}$$

$$\left\|f[t/2]-f^*\right\|_\infty \leq \max\left(\exp\left(-Dt/2\right)\left\|f_0-f^*\right\|_\infty, \max_{l\in[0,t/2]}\left(\left|x(l)-f^*(0)\right|\right)\right) \tag{5.61}$$

Consider the initial-boundary value problem (2.1), (2.2), (2.3) with constant input $D$ and initial condition $\left(f[t/2], S(t/2)\right)$. By virtue of Proposition 1 in [18] (the semigroup property) the solution starting at $t/2$ coincides with the restriction of the original solution of (2.1), (2.2), (2.3), (2.6). So, applying (5.57) and (5.59) with initial time $t/2$ and by (5.60) and (5.61), it follows that

$$\begin{aligned}\left\|f[t]-f^*\right\|_1 &\leq \int_{t/2}^{t} \left|x(l)-f^*(0)\right|dl + \exp\left(-Dt/2\right)\left\|f[t/2]-f^*\right\|_1 \\ &\leq \int_{t/2}^{t} \left|x(l)-f^*(0)\right|dl + \exp\left(-Dt/2\right)\int_0^{t/2} \left|x(l)-f^*(0)\right|dl + \exp\left(-Dt\right)\left\|f_0-f^*\right\|_1\end{aligned} \tag{5.62}$$

$$\left\|f[t]-f^*\right\|_\infty \le \max\left(\exp\left(-Dt/2\right)\left\|f[t/2]-f^*\right\|_\infty, \max_{l\in[t/2,t]}\left(\left|x(l)-f^*(0)\right|\right)\right)$$
$$\le \max\left(\exp\left(-Dt\right)\left\|f_0-f^*\right\|_\infty, \exp\left(-Dt/2\right)\max_{l\in[0,t/2]}\left(\left|x(l)-f^*(0)\right|\right), \max_{l\in[t/2,t]}\left(\left|x(l)-f^*(0)\right|\right)\right) \tag{5.63}$$

Since Lemma 6 establishes the existence of a non-decreasing function $\bar{\rho}:\mathbb{R}_+\to\mathbb{R}_+$ for which (5.31) holds, using definition of $x(t)$ and (5.53) we obtain the following estimate for all $l\in[0,t]$:

$$\begin{aligned}&\left|S(l)-S^*\right|+\left|x(l)-f^*(0)\right|\le \bar{\rho}\left(V\left(f[l],S(l)\right)\right)V^{1/2}\left(f[l],S(l)\right)\\&\le \exp\left(-2\bar{\sigma}\,l\right)\bar{\rho}\left(H\left(f_0,S_0\right)\right)\left(H\left(f_0,S_0\right)\right)^{1/2}\end{aligned} \tag{5.64}$$

Thus, combining estimate (5.62) with (5.64) we get

$$\begin{aligned}&\left\|f[t]-f^*\right\|_1 \le \exp\left(-Dt\right)\left\|f_0-f^*\right\|_1\\&+\bar{\rho}\left(H\left(f_0,S_0\right)\right)\left(H\left(f_0,S_0\right)\right)^{1/2}\int_{t/2}^{t}\exp\left(-2\bar{\sigma}l\right)dl\\&+\exp\left(-Dt/2\right)\bar{\rho}\left(H\left(f_0,S_0\right)\right)\left(H\left(f_0,S_0\right)\right)^{1/2}\int_0^{t/2}\exp\left(-2\bar{\sigma}l\right)dl\\&\le \exp\left(-Dt\right)\left\|f_0-f^*\right\|_1+\bar{\rho}\left(H\left(f_0,S_0\right)\right)\left(H\left(f_0,S_0\right)\right)^{1/2}\frac{\exp\left(-\bar{\sigma}t\right)}{2\bar{\sigma}}\\&+\frac{\exp\left(-Dt/2\right)}{2\bar{\sigma}}\bar{\rho}\left(H\left(f_0,S_0\right)\right)\left(H\left(f_0,S_0\right)\right)^{1/2}\end{aligned} \tag{5.65}$$

Taking into account that $\bar{\sigma}\in\left(0,D/2\right]$, (5.64) and (5.65), we obtain the following estimate:

$$\begin{aligned}&\left(\left|S(t)-S^*\right|+\left\|f[t]-f^*\right\|_1\right)\exp\left(\bar{\sigma}\left(H\left(f_0,S_0\right)\right)t\right)\\&\le\left\|f_0-f^*\right\|_1+\left(H\left(f_0,S_0\right)\right)^{1/2}\bar{\rho}\left(H\left(f_0,S_0\right)\right)\frac{\bar{\sigma}\left(H\left(f_0,S_0\right)\right)+1}{\bar{\sigma}\left(H\left(f_0,S_0\right)\right)}\end{aligned} \tag{5.66}$$

Notice now that (5.64) gives

$$\max_{l\in[0,t/2]}\left(\left|x(l)-f^*(0)\right|\right)\le\bar{\rho}\left(V\left(f_0,S_0\right)\right)\left(V\left(f_0,S_0\right)\right)^{1/2} \tag{5.67}$$

$$\max_{l\in[t/2,t]}\left(\left|x(l)-f^*(0)\right|\right)\le\exp\left(-\bar{\sigma}\,t\right)\bar{\rho}\left(V\left(f_0,S_0\right)\right)\left(V\left(f_0,S_0\right)\right)^{1/2} \tag{5.68}$$

Hence, from (5.67), (5.68), and (5.63), and the fact that $\bar{\sigma} \in (0, D/2]$, the following estimate holds for $t \geq 0$ :

$$\begin{aligned} &\left\| f[t]-f^* \right\|_\infty \leq \max\left( \exp\left(-Dt\right) \left\| f_0 - f^* \right\|_\infty , \max\left( \exp\left(-Dt/2\right), \exp\left(-\bar{\sigma} t\right)\right) G \right) \\ &\leq \exp\left(-\bar{\sigma} t\right) \max\left( \left\| f_0 - f^* \right\|_\infty , G \right) \end{aligned} \tag{5.69}$$

where $G = \bar{\rho}\left(V\left(f_0, S_0\right)\right)\left(V\left(f_0, S_0\right)\right)^{1/2}$. Finally, Lemma 2.4 on page 65 in [17] implies the existence of $\kappa \in K_\infty$ such that

$$s^{1/2} \bar{\rho}\left(s\right) \frac{\bar{\sigma}\left(s\right)+1}{\bar{\sigma}\left(s\right)} \leq \kappa\left(s\right), \text{ for all } s \geq 0 \tag{5.70}$$

Estimates (5.54) and (5.55) follow from (5.66), (5.69) and (5.70). The proof is complete. ⊲

The following lemma provides a useful inequality that is going to be used in the proof of Theorem 1. The inequality is derived with repeated use of the Cauchy-Schwarz and Young inequalities.

**Lemma 8:** *Suppose that assumptions (A), (B) hold. Let $B, \Gamma, M > 0$ and $\sigma \geq 0$ be given constants and define the functional $U : X \to \mathbb{R}$ by means of the following formula for all $\left(f, S\right) \in X$ :*

$$\begin{aligned} &U\left(f,S\right) := -\left( \alpha + \theta \frac{\langle p, v \rangle}{\langle q, r \rangle} \right) \exp(-\varphi)\left(\exp(\varphi)-1\right)^2 - \left(\kappa_1 g\left(S\right)+\delta\right)\left(\exp(\varphi)-1\right)^2 \\ &+\left(\kappa_2-\kappa_1\right)\left(g\left(S\right)-1\right)\left(\exp(\varphi)-1\right) + \frac{\langle \theta p - h, \chi \rangle}{\langle q, r \rangle}\left(\exp(\varphi)-1\right) - \Gamma D g\left(S\right)\left(\exp(w)-1\right)^2 \\ &+\left(\Gamma\kappa_1 - \Gamma D - DM\right)\left(g\left(S\right)-1\right)\left(\exp(w)-1\right) + \frac{B}{2}\left(g\left(S\right)\left(\exp(\varphi)-1\right)+g\left(S\right)-1\right)^2 \\ &+\Gamma\left(\kappa_1 g\left(S\right)+\delta\right)\left(\exp(\varphi)-1\right)\left(\exp(w)-1\right) - DM\frac{\left(g\left(S\right)-1\right)^2}{g\left(S\right)} + \Gamma\frac{\langle h, \chi \rangle}{\langle q, r \rangle}\left(\exp(w)-1\right) \\ &-B\left( \gamma D + \sigma + \left(\kappa_1 g\left(S\right)+\delta\right)\exp(\varphi) + \frac{\langle h, v \rangle}{\langle q, r \rangle} \right)\left\| \rho\chi \right\|_2^2 - B\left\langle \beta\rho^2\chi, \chi \right\rangle \\ &-B\left( \left(\kappa_1 g\left(S\right)+\delta\right)\left(\exp(\varphi)-1\right) + \kappa_1\left(g\left(S\right)-1\right) + \frac{\langle h, \chi \rangle}{\langle q, r \rangle} \right)\left\langle \rho^2 r, \chi \right\rangle \end{aligned} \tag{5.71}$$

*where $v, \chi$ are defined by (3.8). Then for every $\varepsilon, \omega > 0$, $R_1, R_2, R_3 > 0$, $\lambda \in \left[0,1\right]$ the following inequality holds for all $\left(f, S\right) \in X$ :*

$$U(f,S)\le -\left(\alpha-\frac{\omega\lambda}{2\langle q,r\rangle}\right)\exp(-\varphi)\left(\exp(\varphi)-1\right)^2+\left(\kappa_2-\kappa_1\right)\left(g(S)-1\right)\left(\exp(\varphi)-1\right)$$
$$-\left(\kappa_1 g(S)+\delta-Bg^2(S)-B\varepsilon\left(\kappa_1 g(S)+\delta\right)^2-G_2\right)\left(\exp(\varphi)-1\right)^2+G_1\|\rho\chi\|_2^2$$
$$+\left(\Gamma\kappa_1-\Gamma D-DM\right)\left(g(S)-1\right)\left(\exp(w)-1\right)-\Gamma\left(Dg(S)-G_4\right)\left(\exp(w)-1\right)^2$$
$$+\Gamma\delta\frac{R_3}{2}\left(\exp(\varphi)-1\right)^2+\Gamma\kappa_1 g(S)\left(\exp(\varphi)-1\right)\left(\exp(w)-1\right)-B\left(L+\gamma D+\sigma\right)\|\rho\chi\|_2^2$$
$$-B\left(\kappa_1 g(S)+\delta-G_3\right)\exp(\varphi)\|\rho\chi\|_2^2-\left(DM-Bg(S)\left(1+\varepsilon\kappa_1^2\right)\right)\frac{\left(g(S)-1\right)^2}{g(S)} \tag{5.72}$$

*where* $$G_1:=B\left(\frac{\|\rho^{-1}h\|_2}{\langle q,r\rangle}\left(\|\rho r\|_2+\frac{\Gamma}{2BR_2}\right)+\frac{(1-\lambda)\|\rho^{-1}(\theta p-h)\|_2}{2BR_1\langle q,r\rangle}+\frac{\|\rho r\|_2^2}{2\varepsilon}\right),$$

$G_2:=\dfrac{(1-\lambda)R_1}{2\langle q,r\rangle}\|\rho^{-1}(\theta p-h)\|_2$, $G_3:=\dfrac{\lambda\|\rho^{-1}(\theta p-h)\|_2^2}{2\omega B\langle q,r\rangle}$, $G_4:=\dfrac{R_2\|\rho^{-1}h\|_2}{2\langle q,r\rangle}+\dfrac{\delta}{2R_3}$ *and* $L:=\inf\limits_{a\ge 0}\left(\beta(a)\right)$.

**Proof of Lemma 8:** Definition (5.71) and the fact that $h(a)\ge 0$ for all $a\ge 0$ (recall assumption (B)) gives for $L:=\inf\limits_{a\ge 0}\left(\beta(a)\right)$:

$$U(f,S)\le -\alpha\exp(-\varphi)\left(\exp(\varphi)-1\right)^2-\left(\kappa_1 g(S)+\delta\right)\left(\exp(\varphi)-1\right)^2$$
$$+\left(\kappa_2-\kappa_1\right)\left(g(S)-1\right)\left(\exp(\varphi)-1\right)+\frac{\langle\theta p-h,\chi\rangle}{\langle q,r\rangle}\left(\exp(\varphi)-1\right)$$
$$+\Gamma\left(\kappa_1 g(S)+\delta\right)\left(\exp(\varphi)-1\right)\left(\exp(w)-1\right)-\Gamma Dg(S)\left(\exp(w)-1\right)^2$$
$$+\left(\Gamma\kappa_1-\Gamma D-DM\right)\left(g(S)-1\right)\left(\exp(w)-1\right)+\Gamma\frac{\langle h,\chi\rangle}{\langle q,r\rangle}\left(\exp(w)-1\right) \tag{5.73}$$
$$-DM\frac{\left(g(S)-1\right)^2}{g(S)}+\frac{B}{2}\left(g(S)\left(\exp(\varphi)-1\right)+g(S)-1\right)^2$$
$$-B\left(L+\gamma D+\sigma+\left(\kappa_1 g(S)+\delta\right)\exp(\varphi)\right)\|\rho\chi\|_2^2-B\frac{\langle h,\chi\rangle}{\langle q,r\rangle}\left\langle\rho^2 r,\chi\right\rangle$$
$$-B\left(\left(\kappa_1 g(S)+\delta\right)\left(\exp(\varphi)-1\right)+\kappa_1\left(g(S)-1\right)\right)\left\langle\rho^2 r,\chi\right\rangle$$

Using the Cauchy-Schwarz and Young inequalities we conclude that the following inequalities hold for all $\varepsilon>0$:

$$\left|\langle h,\chi\rangle\langle\rho^2 r,\chi\rangle\right|\le\left\|\rho^{-1}h\right\|_2\|\rho r\|_2\|\rho\chi\|_2^2$$

$$\frac{1}{2}\left(g(S)\left(\exp(\varphi)-1\right)+g(S)-1\right)^2\le g^2(S)\left(\exp(\varphi)-1\right)^2+\left(g(S)-1\right)^2 \tag{5.74}$$

$$\left|\left(\left(\kappa_1 g(S)+\delta\right)\left(\exp(\varphi)-1\right)+\kappa_1\left(g(S)-1\right)\right)\langle\rho r,\chi\rangle\right|$$
$$\le\frac{\varepsilon}{2}\left(\left(\kappa_1 g(S)+\delta\right)\left(\exp(\varphi)-1\right)+\kappa_1\left(g(S)-1\right)\right)^2+\frac{1}{2\varepsilon}\|\rho r\|_2^2\|\rho\chi\|_2^2$$

Inequalities (5.74) in conjunction with (5.73) give for any $\lambda\in[0,1]$:

$$\begin{aligned}
U(f,S)\le&-\alpha\exp(-\varphi)\left(\exp(\varphi)-1\right)^2-\left(\kappa_1 g(S)+\delta-Bg^2(S)\right)\left(\exp(\varphi)-1\right)^2\\
&+\left(\kappa_2-\kappa_1\right)\left(g(S)-1\right)\left(\exp(\varphi)-1\right)+(1-\lambda)\frac{\langle\theta p-h,\chi\rangle}{\langle q,r\rangle}\left(\exp(\varphi)-1\right)\\
&+\lambda\frac{\langle\theta p-h,\chi\rangle}{\langle q,r\rangle}\left(\exp(\varphi)-1\right)+\Gamma\left(\kappa_1 g(S)+\delta\right)\left(\exp(\varphi)-1\right)\left(\exp(w)-1\right)\\
&-\Gamma Dg(S)\left(\exp(w)-1\right)^2+\left(\Gamma\kappa_1-\Gamma D-DM\right)\left(g(S)-1\right)\left(\exp(w)-1\right)+\Gamma\frac{\langle h,\chi\rangle}{\langle q,r\rangle}\left(\exp(w)-1\right)\\
&-\left(DM-Bg(S)\right)\frac{\left(g(S)-1\right)^2}{g(S)}+B\frac{\varepsilon}{2}\left(\left(\kappa_1 g(S)+\delta\right)\left(\exp(\varphi)-1\right)+\kappa_1\left(g(S)-1\right)\right)^2\\
&-B\left(L+\gamma D+\sigma+\left(\kappa_1 g(S)+\delta\right)\exp(\varphi)-\frac{\left\|\rho^{-1}h\right\|_2}{\langle q,r\rangle}\|\rho r\|_2-\frac{1}{2\varepsilon}\|\rho r\|_2^2\right)\|\rho\chi\|_2^2
\end{aligned} \tag{5.75}$$

By virtue of Young and Cauchy-Schwarz inequalities, we get for arbitrary $\omega>0$:

$$\frac{1}{2}\left(\left(\kappa_1 g(S)+\delta\right)\left(\exp(\varphi)-1\right)+\kappa_1\left(g(S)-1\right)\right)^2$$
$$\le\left(\kappa_1 g(S)+\delta\right)^2\left(\exp(\varphi)-1\right)^2+\kappa_1^2\left(g(S)-1\right)^2$$

$$\langle\theta p-h,\chi\rangle\left(\exp(\varphi)-1\right)\le\frac{\omega}{2}\exp(-\varphi)\left(\exp(\varphi)-1\right)^2+\frac{1}{2\omega}\exp(\varphi)\left\|\rho^{-1}\left(\theta p-h\right)\right\|_2^2\|\rho\chi\|_2^2 \tag{5.76}$$

Combining, (5.76) and (5.75) we obtain

$$U(f,S) \le -\left(\alpha - \frac{\omega\lambda}{2\langle q,r\rangle}\right)\exp(-\varphi)\left(\exp(\varphi)-1\right)^2 - \Gamma D g(S)\left(\exp(w)-1\right)^2$$

$$-\left(\kappa_1 g(S)+\delta - Bg^2(S) - B\varepsilon\left(\kappa_1 g(S)+\delta\right)^2\right)\left(\exp(\varphi)-1\right)^2 + \Gamma\frac{\langle h,\chi\rangle}{\langle q,r\rangle}\left(\exp(w)-1\right)$$

$$+\Gamma\left(\kappa_1 g(S)+\delta\right)\left(\exp(\varphi)-1\right)\left(\exp(w)-1\right) - \left(DM - Bg(S)\left(1+\varepsilon\kappa_1^2\right)\right)\frac{\left(g(S)-1\right)^2}{g(S)}$$

$$+\left(\Gamma\kappa_1 - \Gamma D - DM\right)\left(g(S)-1\right)\left(\exp(w)-1\right) + (1-\lambda)\frac{\langle \theta p - h,\chi\rangle}{\langle q,r\rangle}\left(\exp(\varphi)-1\right)$$

$$+\left(\kappa_2-\kappa_1\right)\left(g(S)-1\right)\left(\exp(\varphi)-1\right) - B\left(\kappa_1 g(S)+\delta - \frac{\lambda\left\|\rho^{-1/2}\left(\theta p - h\right)\right\|_2^2}{2\omega B\langle q,r\rangle}\right)\exp(\varphi)\left\|\rho\chi\right\|_2^2$$

$$-B\left(L+\gamma D+\sigma - \frac{\left\|\rho^{-1}h\right\|_2}{\langle q,r\rangle}\left\|\rho r\right\|_2 - \frac{1}{2\varepsilon}\left\|\rho r\right\|_2^2\right)\left\|\rho\chi\right\|_2^2 \tag{5.77}$$

Finally, for arbitrary $R_1,R_2,R_3>0$ we get (using Young and Cauchy-Schwarz inequalities) the following inequalities

$$\langle \theta p - h,\chi\rangle\left(\exp(\varphi)-1\right) \le \frac{R_1}{2}\left\|\rho^{-1}\left(\theta p-h\right)\right\|_2\left(\exp(\varphi)-1\right)^2 + \frac{1}{2R_1}\left\|\rho^{-1}\left(\theta p-h\right)\right\|_2\left\|\rho\chi\right\|_2^2$$

$$\langle h,\chi\rangle\left(\exp(w)-1\right) \le \frac{R_2}{2}\left\|\rho^{-1}h\right\|_2\left(\exp(w)-1\right)^2 + \frac{1}{2R_2}\left\|\rho^{-1}h\right\|_2\left\|\rho\chi\right\|_2^2 \tag{5.78}$$

$$\left(\exp(\varphi)-1\right)\left(\exp(w)-1\right) \le \frac{R_3}{2}\left(\exp(\varphi)-1\right)^2 + \frac{1}{2R_3}\left(\exp(w)-1\right)^2$$

which, when applied to (5.77) give inequality (5.72). The proof is complete. ◁

We are now ready to provide the proof of Theorem 1.

**Proof of Theorem 1:** We first focus on the case where the initial condition is in the set $\Omega$ defined by (3.6).

Let arbitrary $\left(f_0,S_0\right)\in\Omega$ be given. By Lemma 2, the corresponding weak solution of (2.1), (2.2), (2.3), (2.6) exists and satisfies $\left(f[t],S(t)\right)\in\Omega$ for all $t\ge 0$. Using (5.26),

(5.22), (3.12), (3.13), (3.10) and definition (5.71) we obtain the following equation for all $t \geq 0$:

$$\frac{d}{dt}\left(V(f[t],S(t))\right) = U(f[t],S(t)) \tag{5.79}$$

It follows from Lemma 8 and equation (5.79) that for every $\varepsilon, \omega > 0$, $R_1, R_2, R_3 > 0$, $\lambda \in [0,1]$ the following inequality holds for all $t \geq 0$:

$$\begin{aligned}
&\frac{d}{dt}\left(V(f[t],S(t))\right) \leq -\left(\alpha - \frac{\omega\lambda}{2\langle q,r\rangle}\right)\exp(-\varphi(t))\left(\exp(\varphi(t))-1\right)^2 + G_1\left\|\rho\chi[t]\right\|_2^2 \\
&+\left(\kappa_2-\kappa_1\right)\left(g\left(S(t)\right)-1\right)\left(\exp(\varphi(t))-1\right) - B\left(\kappa_1 g\left(S(t)\right)+\delta-G_3\right)\exp(\varphi(t))\left\|\rho\chi[t]\right\|_2^2 \\
&-\left(\kappa_1 g\left(S(t)\right)+\delta-Bg^2\left(S(t)\right)-B\varepsilon\left(\kappa_1 g\left(S(t)\right)+\delta\right)^2-G_2\right)\left(\exp(\varphi(t))-1\right)^2 \\
&+\left(\Gamma\kappa_1-\Gamma D-DM\right)\left(g\left(S(t)\right)-1\right)\left(\exp(\varphi(t))-1\right)+\Gamma\delta\frac{R_3}{2}\left(\exp(\varphi(t))-1\right)^2 \\
&+\Gamma\kappa_1 g\left(S(t)\right)\left(\exp(\varphi(t))-1\right)\left(\exp(w(t))-1\right)-\Gamma\left(Dg\left(S(t)\right)-G_4\right)\left(\exp(w(t))-1\right)^2 \\
&-\left(DM-Bg\left(S(t)\right)\left(1+\varepsilon\kappa_1^2\right)\right)\frac{\left(g\left(S(t)\right)-1\right)^2}{g\left(S(t)\right)} - B\left(L+\gamma D+\sigma\right)\left\|\rho\chi[t]\right\|_2^2
\end{aligned} \tag{5.80}$$

Notice now that since $(f[t],S(t)) \in \Omega$, definition (3.6) gives for all $t \geq 0$:

$$0 < g\left(\underline{S}\right) \leq g\left(S(t)\right) \leq g\left(S_{in}\right) \tag{5.81}$$

Selecting $R_1, R_2, R_3 > 0$, $\omega > 0$, $\varepsilon > 0$, $\lambda \in [0,1]$, $B, \Gamma, M > 0$ so that conditions (3.20), (3.21), (3.22), and (3.23) are valid, together with (5.80), (5.81) and definition (3.24), we get for all $t \geq 0$:

$$\frac{d}{dt}\left(V(f[t],S(t))\right) \leq -z^T(t)Pz(t) - c\left\|\rho\chi[t]\right\|_2^2 \leq 0 \tag{5.82}$$

where $z^T(t)$ is the transpose of

$$z(t) = \begin{pmatrix} \sqrt{g(S(t))}\left(e^{\varphi(t)}-1\right) \\ \sqrt{g(S(t))}\left(e^{w(t)}-1\right) \\ (g(S(t))-1)/\sqrt{g(S(t))} \end{pmatrix} \tag{5.83}$$

and

$$\begin{aligned}
&c := L+\gamma D+\sigma \\
&-\left(\frac{\left\|\rho^{-1}h\right\|_2}{\langle q,r\rangle}\left(\left\|\rho r\right\|_2+\frac{\Gamma}{2BR_2}\right)+\frac{(1-\lambda)\left\|\rho^{-1}\left(\theta p-h\right)\right\|_2}{2BR_1\langle q,r\rangle}+\frac{\left\|\rho r\right\|_2^2}{2\varepsilon}\right) > 0
\end{aligned}$$

Define:

$$k_0 = \min\left( \lambda_{\min}\left(P\right) \min\left( g(\underline{S}), \frac{1}{g(S_{in})} \right), c \right) > 0 \tag{5.84}$$

where $\lambda_{\min}\left(P\right) > 0$ is the smallest eigenvalue of the positive definite matrix $P$ in (3.24). Then, (5.82) and Lemma 6 give

$$\frac{d}{dt}\left(V(f[t],S(t))\right) \leq -\frac{k_0 V\left(f[t],S(t)\right)}{\bar{\gamma}\left(V\left(f[t],S(t)\right)\right)} \tag{5.85}$$

for an appropriate non-decreasing function $\bar{\gamma} : \mathbb{R}_+ \to (0,+\infty)$. Notice that inequality (5.82) implies that $V(f[t],S(t)) \leq V(f[0],S(0)) = V(f_0,S_0)$ for all $t \geq 0$. Thus, (5.85) and the fact that $\bar{\gamma} : \mathbb{R}_+ \to (0,+\infty)$ is a non-decreasing function imply the following differential inequality for all $t \geq 0$:

$$\frac{d}{dt}\left(V(f[t],S(t))\right) \leq -\frac{k_0}{\bar{\gamma}\left(V\left(f_0,S_0\right)\right)} V\left(f[t],S(t)\right) \tag{5.86}$$

The differential inequality (5.86) gives for all $t \geq 0$:

$$V\left(f[t],S(t)\right) \leq \exp\left( -\frac{k_0\, t}{\bar{\gamma}\left(V\left(f_0,S_0\right)\right)} \right) V\left(f_0,S_0\right) \tag{5.87}$$

We next focus on the case where the initial condition is not in the set $\Omega$ defined by (3.6).

Let arbitrary $\left(f_0,S_0\right) \in X \setminus \Omega$ be given. Lemma 2 and definition (3.6) guarantee the existence of an increasing continuous function $T : \mathbb{R}_+ \to \mathbb{R}_+$ such that the weak solution $\left(f[t],S(t)\right)$ of the initial-boundary value problem (2.1), (2.2), (2.3), (2.6) satisfies $\left(f[t],S(t)\right) \in \Omega$ for all $t \geq T\left(RS_0 + \|f_0\|_1\right)$. It follows from the semigroup property and (5.87) that the following estimate holds for all $t \geq T := T\left(RS_0 + \|f_0\|_1\right)$:

$$V\left(f[t],S(t)\right) \leq \exp\left( -\frac{k_0\left(t-T\right)}{\bar{\gamma}\left(V\left(f[T],S(T)\right)\right)} \right) V\left(f[T],S(T)\right) \tag{5.88}$$

It follows from (2.5), (3.7), (2.17), (2.18), (2.19), (5.19), (5.20), (5.21) that the following estimates hold for all $t \in \left[0,T\right]$:

$$\left|\zeta(t) - 1\right| \leq \frac{\|q\|_\infty \exp\left(\bar{R}\,T\right)}{\left\langle q, f^* \right\rangle} d\left(\left(f_0,S_0\right),\left(f^*,S^*\right)\right) \tag{5.89}$$

$$|\xi(t)-1| \le \frac{\|k\|_\infty \exp\left(\bar{R}T\right)}{\left\langle k, f^*\right\rangle} d\left(\left(f_0, S_0\right), \left(f^*, S^*\right)\right) \tag{5.90}$$

$$d\left(\left(f[t], S(t)\right), \left(f^*, S^*\right)\right) \le \exp\left(\bar{R}T\right) d\left(\left(f_0, S_0\right), \left(f^*, S^*\right)\right) \tag{5.91}$$

$$\left\|f[t]-f^*\right\|_\infty \le \max\left(\left\|f_0 - f^*\right\|_\infty, C\exp\left(\bar{R}T\right) d\left(\left(f_0, S_0\right), \left(f^*, S^*\right)\right)\right) \tag{5.92}$$

$$S(t) \le S_{in} - \left(S_{in} - S_0\right)\exp\left(-DT\right) \tag{5.93}$$

$$\xi(t) \ge \exp\left(D\left(b-1\right)T\right)\frac{\left\langle k, f_0\right\rangle}{\left\langle k, f^*\right\rangle} \tag{5.94}$$

$$\zeta(t) \ge \exp\left(D\left(\gamma-1\right)T\right)\frac{\left\langle q, f_0\right\rangle}{\left\langle q, f^*\right\rangle} \tag{5.95}$$

$$S(t) \ge \frac{DS_0}{D + L_\mu \|q\|_\infty \left(RS_{in} + \left(RS_0 + \|f_0\|_1 - RS_{in}\right)^+\right)} \tag{5.96}$$

Definitions (3.9), inequality (3.1) and estimates (5.89), (5.90), (5.93), (5.94), (5.95), (5.96) imply that that the following estimates hold for all $t \in \left[0, T\right]$:

$$\varphi(t) \ge \left(D\left(b-1\right) - \bar{R}\right)T + \ln\left(\frac{\left\langle q, f^*\right\rangle\left\langle k, f_0\right\rangle}{\left\langle k, f^*\right\rangle\left(\left\langle q, f^*\right\rangle + \|q\|_\infty d\left(\left(f_0, S_0\right), \left(f^*, S^*\right)\right)\right)}\right) \tag{5.97}$$

$$\varphi(t) \le \ln\left(R\frac{\left\langle q, f^*\right\rangle}{\left\langle k, f^*\right\rangle}\right) \tag{5.98}$$

$$w(t) \ge D\left(\gamma-1\right)T + \ln\left(\frac{\left(S_{in} - S^*\right)\left\langle q, f_0\right\rangle}{S_{in}\left\langle q, f^*\right\rangle}\right) \tag{5.99}$$

$$w(t) \le \left(\bar{R} + D\right)T + \ln\left(\frac{S_{in} - S^*}{S_{in} - S_0}\left(1 + \frac{\|q\|_\infty}{\left\langle q, f^*\right\rangle} d\left(\left(f_0, S_0\right), \left(f^*, S^*\right)\right)\right)\right) \tag{5.100}$$

It also follows from (3.8), (3.16) with $\sigma \ge 0$ and by using Young and Cauchy-Schwarz inequalities that the following inequalities hold for all $t \ge 0$

$$
\begin{aligned}
\|\rho\chi[t]\|_2^2 &= \frac{1}{\zeta^2(t)\left(f^*(0)\right)^2}\left\|\rho\left(f[t]-f^*\right)\right\|_2^2 \\
&+\frac{\left(\zeta(t)-1\right)^2}{\zeta^2(t)}\|\rho r\|_2^2-\frac{2\left(\zeta(t)-1\right)}{\zeta^2(t)f^*(0)}\left\langle \rho r,\rho\left(f[t]-f^*\right)\right\rangle \\
&\le \frac{2}{\zeta^2(t)\left(f^*(0)\right)^2}\left\|\rho\left(f[t]-f^*\right)\right\|_2^2+2\frac{\left(\zeta(t)-1\right)^2}{\zeta^2(t)}\|\rho r\|_2^2 \\
&\le \frac{2\left\|f[t]-f^*\right\|_\infty}{\zeta^2(t)\left(f^*(0)\right)^2}\left\|f[t]-f^*\right\|_1+2\frac{\left(\zeta(t)-1\right)^2}{\zeta^2(t)}\|\rho r\|_2^2
\end{aligned}
\tag{5.101}
$$

Combining (2.5), (3.15), (5.101), (5.97), (5.98), (5.99), (5.100), (5.93), (5.96), (5.91), (5.92), (5.89) and (5.95), we can conclude that there exists a functional $\bar{V}: X \to \mathbb{R}_+$ such that the following estimate holds for all $t \in [0,T]$:

$$V\left(f[t],S(t)\right) \le \bar{V}\left(f_0,S_0\right) \tag{5.102}$$

Using the same arguments as in the proof of Lemma 6, we can guarantee that the functional $\bar{V}: X \to \mathbb{R}_+$ satisfies the following property for all $r \ge 0$:

$$\sup\left\{ \bar{V}\left(f,S\right):\left(f,S\right)\in X, \left\|f-f^*\right\|_\infty + V\left(f,S\right)+\left\|f-f^*\right\|_1 \le r \right\} < +\infty \tag{5.103}$$

Using (5.102) and (5.88) we conclude that there exists a functional $\tilde{V}: X \to \mathbb{R}_+$ such that the following estimates hold for all $t \ge 0$:

$$V\left(f[t],S(t)\right) \le \exp\left(-\frac{k_0}{\bar{\gamma}\left(\tilde{V}\left(f_0,S_0\right)\right)}t\right)\tilde{V}\left(f_0,S_0\right) \tag{5.104}$$

Combining estimates (5.87) (for the case $\left(f_0,S_0\right)\in\Omega$), (5.102) and (5.88) (for the case $\left(f_0,S_0\right)\in X\setminus\Omega$), we conclude that there exist a non-increasing function $\bar{\sigma}:\mathbb{R}_+ \to (0,+\infty)$ and a functional $H: X \to \mathbb{R}_+$ with $H = V$ on $\Omega$ that satisfies for all $r \ge 0$

$$\sup\left\{ H\left(f,S\right):\left(f,S\right)\in X, \left\|f-f^*\right\|_\infty + V\left(f,S\right)+\left\|f-f^*\right\|_1 \le r \right\} < +\infty \tag{5.105}$$

and such that for every $\left(f_0,S_0\right)\in X$ the weak solution of the initial-boundary value problem (2.1), (2.2), (2.3), (2.6) satisfies estimate (5.53) for all $t \ge 0$ with $\bar{\sigma} := \bar{\sigma}\left(H\left(f_0,S_0\right)\right)$. It follows from Lemma 7 and definition (3.18) that there exists a function $\kappa \in K_\infty$ such that for every $\left(f_0,S_0\right)\in X$ the weak solution of the initial-boundary value problem (2.1), (2.2), (2.3), (2.6) satisfies the following estimate for all $t \ge 0$:

$$\Psi\left(f[t],S(t)\right)\leq\exp\left(-\bar{\sigma}\left(H\left(f_0,S_0\right)\right)t\right)\left(\left\|f_0-f^*\right\|_1+\left\|f_0-f^*\right\|_\infty+\kappa\left(H\left(f_0,S_0\right)\right)\right) \tag{5.106}$$

It is straightforward to show that $\left(f^*,S^*\right)\in\Omega$. Indeed, exploiting (2.10) and (2.11) we get $RS^*+\left\|f^*\right\|_1=RS^*+\frac{D\langle k,r\rangle\left(S_{in}-S^*\right)}{\langle q,r\rangle}\|r\|_1$. Definition (2.12) and the fact that $\beta$ is non-negative implies that $\|r\|_1=\int_0^{+\infty}\exp\left(-Da-\int_0^a\beta(s)ds\right)da\leq\frac{1}{D}$. Combining the previous relations with (3.1) (which implies that $\frac{\langle k,r\rangle}{\langle q,r\rangle}\leq R$) gives

$$RS^*+\left\|f^*\right\|_1\leq RS_{in}<F \tag{5.107}$$

Similarly, exploiting (2.11) (which gives $S^*=S_{in}-\frac{\mu\left(S^*\right)\langle q,r\rangle}{D\|r\|_1}\left\|f^*\right\|_1$), inequality (5.107) (which implies that $\left\|f^*\right\|_1\leq F$) and the inequality $\langle q,r\rangle\leq\|q\|_\infty\|r\|_1$ we obtain the inequality $S^*\geq S_{in}-\frac{\mu\left(S^*\right)\|q\|_\infty}{D}F$. Definition (3.5) and the fact that $\mu(0)=0$ gives $\mu\left(S^*\right)\leq L_\mu S^*$. Combining the previous relations gives

$$S^*\geq\frac{DS_{in}}{D+L_\mu F\|q\|_\infty}=2\underline{S} \tag{5.108}$$

Therefore, definition (3.6) and (5.107), (5.108) imply that $\left(f^*,S^*\right)\in\Omega$. Moreover, inequalities (5.107), (5.108) and definitions (3.6), (3.18) guarantee that the following implication holds:

$$\left(f,S\right)\in X,\Psi\left(f,S\right)\leq\min\left(\underline{S},\frac{F-RS_{in}}{1+R}\right)\Rightarrow\left(f,S\right)\in\Omega \tag{5.109}$$

We define for all $r\geq 0$:

$$\bar{b}(r)$$
$$:=\sup\left\{\left\|f-f^*\right\|_\infty+\left\|f-f^*\right\|_1+H\left(f,S\right)+\kappa\left(H\left(f,S\right)\right):\left(f,S\right)\in X,\Psi\left(f,S\right)\leq r\right\} \tag{5.110}$$

Inequality (5.105) shows that the function $\bar{b}$ is well-defined on $\mathbb{R}_+$, i.e., $0\leq\bar{b}(r)<+\infty$ for all $r\geq 0$. Moreover, definition (5.110) guarantees that the function $\bar{b}$ is non-decreasing on $\mathbb{R}_+$. Furthermore, the fact that $H = V$ on $\Omega$, implication (5.109) and definition (3.18) guarantee the following implication:

$$\Psi(f,S) \le r \le \min\left( \underline{S}, \frac{F - RS_{in}}{1+R} \right) \Rightarrow 0 \le \bar{b}(r) \le r + \kappa(r) \tag{5.111}$$

Implication (5.111) guarantees that $0 = \bar{b}(0) = \limsup_{r \to 0^+} \left( \bar{b}(r) \right)$. Lemma 2.4 on page 65 in [17] implies the existence of $C \in K_\infty$ such that

$$\bar{b}(r) \le C(r), \text{ for all } r \ge 0 \tag{5.112}$$

Estimate (3.25) is a direct consequence of estimate (5.106), definition (5.110) and inequality (5.112). The proof is complete. $\triangleleft$

**Proof of Proposition 1:** Suppose that condition (4.9) is satisfied. We will verify that all conditions of Theorem 1 are satisfied for system (4.1), (4.2), (4.3). Then, by Theorem 1, estimate (3.25) holds, and the equilibrium $(f^*, S^*)$ is globally asymptotically stable. Using definitions (3.10) and (3.19) we get

$$\theta = \frac{D + L + \tilde{k}}{Y(D+L)} \tag{5.113}$$

$$\begin{aligned} \kappa_1 &= D + L, \\ \kappa_2 &= D + L + \tilde{k} \end{aligned} \tag{5.114}$$

Selecting

$$\begin{aligned} &\gamma = -\frac{L}{D}, b = -\frac{k+L}{D}, \\ &\delta = \alpha = 0, \\ &h(a) \equiv 0, p(a) \equiv 0 \end{aligned} \tag{5.115}$$

assumption (B) is satisfied while assumption (A) can be verified for any $R \ge Y$. With the above selection and by also setting $\lambda = 0$, conditions (3.20), (3.21), and (3.23) hold. Using (5.115), inequality (3.22) is satisfied for

$$\varepsilon > \frac{1}{4\sigma(\sigma + D + L)} \text{ and } \sigma > 0 \tag{5.116}$$

Next, using (5.114) and (5.115), the matrix $P$ in (3.24) is given by

$$P = \begin{bmatrix} L + D - Bg(S_{in})(1+\varepsilon\kappa_1^2) & -\Gamma(L+D)/2 & -\tilde{k}/2 \\ -\Gamma(L+D)/2 & \Gamma D & (DM - \Gamma L)/2 \\ -\tilde{k}/2 & (DM - \Gamma L)/2 & DM - Bg(S_{in})(1+\varepsilon\kappa_1^2) \end{bmatrix} \tag{5.117}$$

To apply Theorem 1, it suffices to show the existence of constants $B,\Gamma,M>0$ such that the matrix $P$ in (5.117) is positive definite.

Notice first that, for sufficiently small $B>0$, the matrix $P$ in (5.117) is positive definite if the following matrix is positive definite

$$P_0=\begin{bmatrix} L+D & -\Gamma(L+D)/2 & -\tilde{k}/2 \\ -\Gamma(L+D)/2 & \Gamma D & \left(DM-\Gamma L\right)/2 \\ -\tilde{k}/2 & \left(DM-\Gamma L\right)/2 & DM \end{bmatrix} \tag{5.118}$$

Indeed, if $P_0$ is positive definite, its smallest eigenvalue $\lambda_{\min}\left(P_0\right)$ is positive and therefore $P$ remains positive definite if $0<B<\lambda_{\min}(P_0)/\left(g(S_{in})(1+\varepsilon\kappa_1^2\right)$. Thus, it suffices to verify positive definiteness of the matrix $P_0$ in (5.118) by appropriately selecting the constants $M$, $\Gamma$. Notice that since $L,D>0$, by Sylvester's criterion, $P_0$ is positive definite if and only if

$$\begin{aligned} &\Gamma(L+D)<4D \\ &J(M)>\tilde{k}^2\Gamma D \end{aligned} \tag{5.119}$$

where

$$J(M)=(L+D)\left(-D^2M^2+D\Gamma\left(4D+2L-\Gamma(L+D)+\tilde{k}\right)M-\Gamma^2L(L+\tilde{k})\right)$$

Notice now that $J(M)$ above is a strictly concave quadratic polynomial in $M$ and thus, it attains a unique maximum when

$$M=\frac{4\Gamma D-\Gamma^2(L+D)+2\Gamma L+\tilde{k}\Gamma}{2D} \tag{5.120}$$

with $M>0$ when $\Gamma(L+D)<4D$. Thus, from (5.119) and (5.120), positive definiteness of $P_0$ reduces to the conditions

$$\begin{aligned} &\Gamma(L+D)<4D \\ &\tilde{E}(\Gamma):=\Gamma L\left(4D-\Gamma(L+D)\right)+\frac{\Gamma\left(4D-\Gamma(L+D)+\tilde{k}\right)^2}{4}-\frac{\tilde{k}^2D}{L+D} \\ &=\frac{\left(\Gamma(L+D)-4D\right)\left((L+D)^2\Gamma^2-2(L+D)(2(L+D)+\tilde{k})\Gamma+\tilde{k}^2\right)}{4(L+D)}>0 \end{aligned} \tag{5.121}$$

Notice now that $\tilde{E}\left(\frac{4D}{L+D}\right)=0$ and that for $\Gamma\in\left(0,\frac{4D}{L+D}\right)$ we have $\Gamma(L+D)<4D$, thus, from (5.121) it follows that $\tilde{E}(\Gamma)>0$ if and only if

$$\Delta(\Gamma):=(L+D)^2\Gamma^2-2(L+D)(2(L+D)+\tilde{k})\Gamma+\tilde{k}^2<0 \tag{5.122}$$

for some $\Gamma\in\left(0,\frac{4D}{L+D}\right)$. Inequality (5.122), shows that $\Delta(\Gamma)<0$ for all $\Gamma\in\left(\Gamma_1,\Gamma_2\right)$ where

$$\Gamma_1 = \left( \sqrt{\frac{\tilde{k}}{L+D}+1} - 1 \right)^2 , \qquad \Gamma_2 = \left( \sqrt{\frac{\tilde{k}}{L+D}+1} + 1 \right)^2 \tag{5.123}$$

Since $\tilde{k} \geq 0$ and $L+D > D > 0$, it follows that $\Gamma_2 \geq 4\sqrt{\frac{\tilde{k}}{L+D}+1} > \frac{4D}{L+D}$. Thus, (5.121) is satisfied for some $\Gamma \in \left(0, \frac{4D}{D+L}\right)$ if and only if $\Gamma_1 < \frac{4D}{D+L}$. Let now $\tilde{g}(D) := \frac{4D}{L+D} - \left( \sqrt{\frac{\tilde{k}}{L+D}+1} - 1 \right)^2$, $D > 0$ and notice that $\tilde{g}$ is strictly increasing in $D$ with $\tilde{g}\left( \frac{\tilde{k}^2}{8(2L+\tilde{k})} \right) = 0$. Thus, $\tilde{g}(D) > 0$ for any $D$ satisfying (4.9). The latter implies that there exists $\Gamma \in \left(0, \frac{4D}{D+L}\right)$ such that $\Gamma_1 < \frac{4D}{D+L}$. Consequently, we have found $B, \Gamma, M > 0$ such that the matrix $P$ in (5.117) is positive definite. Hence, under inequality (4.9), all conditions of Theorem 1 are satisfied. The proof is complete. $\triangleleft$

## 6. Conclusions

In this paper, we have derived global $KL$ stability estimates for an age-structured chemostat model coupled with substrate dynamics under a constant dilution rate. Under appropriate structural assumptions linking the birth and substrate consumption rates, we established the existence of a trapping region that attracts all trajectories in finite time. We then constructed an appropriate Lyapunov functional, and derived sufficient conditions ensuring global exponential decay along trajectories.

Beyond equilibrium analysis, chemostat models can also be viewed from a control-theoretic perspective, since the dilution rate can be manipulated in practice and act as a control input (see for instance [2], [6], [13], [9], [16], [20], [21], [19], [28]). Stabilization of chemostat systems by means of dilution control has been extensively studied both in the classical finite-dimensional setting and for age-structured models (under strong structural assumptions or in the absence of the substrate equations), see for instance [4], [8], [14], [15], [22], [34], [35]. Although the present work focuses on constant dilution rates, the formulation of the model as a control system in [18] and the construction of a Lyapunov functional provide a natural foundation for future feedback stabilization results in the age-structured setting.

## Appendix

**Proof of Lemma 5:** In what follows $\nabla^2 W(x)$ denotes the Hessian matrix at $x \in O \subseteq \mathbb{R}^n$ of the function $W \in C^2(O)$. Since $W \in C^2(O)$ is a positive definite function with $\nabla^2 W(0)$ being positive definite, there exist $c, \delta > 0$ such that the following implication holds:

$$x \in O, |x| \le \delta \Rightarrow W(x) \ge c|x|^2 \tag{A.1}$$

Indeed, since $\nabla^2 W(0)$ is positive definite, there exists $c > 0$ such that $\xi' \nabla^2 W(0) \xi \ge 4c|\xi|^2$ for all $\xi \in \mathbb{R}^n$. Since $W \in C^2(D)$ there exists $\delta > 0$ such that $\left|\nabla^2 W(z) - \nabla^2 W(0)\right| \le 2c$ for all $z \in O$ with $|z| \le \delta$. The fact that $W \in C^2(O)$ is a positive definite function implies that $W(0) = 0$, $\nabla W(0) = 0$ and $W(x) = \int_0^1 \int_0^1 \lambda x' \nabla^2 W(\lambda \mu x) x \, d\mu d\lambda$ for all $x \in D$ with $|x| \le \delta$. Clearly, we obtain $W(x) = \int_0^1 \int_0^1 \lambda x' \nabla^2 W(0) x \, d\mu d\lambda + \int_0^1 \int_0^1 \lambda x' \left(\nabla^2 W(\lambda \mu x) - \nabla^2 W(0)\right) x \, d\mu d\lambda$, which gives $W(x) \ge 2c|x|^2 - |x|^2 \int_0^1 \int_0^1 \lambda \left|\nabla^2 W(\lambda \mu x) - \nabla^2 W(0)\right| d\mu d\lambda$. Implication (A.1) is a straightforward consequence of the above inequalities.

We define $a := \frac{1}{2}\inf\left\{ W(x) : x \in O, |x| \geq \delta \right\}$ (since $O \subseteq \mathbb{R}^n$ is open with $0 \in D$, it follows that the set $\left\{ x \in O : |x| \geq \delta \right\}$ is non-empty for sufficiently small $\delta > 0$). Notice that the definition of $a$ implies that $a \geq 0$ (since $W$ is non-negative). We prove with contradiction that $a > 0$. If $a$ were equal to zero then we would obtain a sequence $\left\{ x_n \in O : n = 1,2,... \right\}$ with $|x_n| \geq \delta$ for $n = 1,2,...$ such that $\lim\left(W\left(x_n\right)\right) = 0$. It follows that the sequence $\left\{ W\left(x_n\right) : n = 1,2,... \right\}$ is bounded and there exists $r \geq 0$ such that $W\left(x_n\right) \leq r$ for $n = 1,2,...$. Consequently, the sequence $\left\{ x_n \in D : n = 1,2,... \right\}$ is in the compact set $\left\{ x \in O : W(x) \leq r \right\}$. Therefore, the Bolzano–Weierstrass theorem implies the existence of a convergent subsequence $\left\{ y_k \in O : k = 1,2,... \right\} \subseteq \left\{ x_n \in O : n = 1,2,... \right\}$ with $|y_k| \geq \delta$ for $k = 1,2,...$, $\lim\left(W\left(y_k\right)\right) = 0$ and $\lim\left(y_k\right) = y^* \in O$. This implies $\left|y^*\right| \geq \delta$ and $W\left(y^*\right) = 0$; a contradiction.

We next show that the following implication holds:

$$x \in O, W(x) \leq a \Rightarrow W(x) \geq c|x|^2 \tag{A.2}$$

Indeed, let arbitrary $x \in O$ with $W(x) \leq a$ be given. If $|x| > \delta$ then the definition of $a$ would imply that $W(x) \geq 2a$, contradicting the fact that $W(x) \leq a$. Thus, we must have $|x| \leq \delta$ and implication (A.1) shows that $W(x) \geq c|x|^2$.

Define the non-decreasing function $q : \mathbb{R}_+ \to \mathbb{R}_+$ by means of the formula $q\left(r\right) := \max\left\{ |\xi|^2 : \xi \in O, W(\xi) \leq r \right\}$ for $r \geq 0$. Implication (A.2) shows that $q\left(r\right) \leq c^{-1} r$ for all $r \in \left[0, a\right]$. Moreover, we have $|x|^2 \leq q(W(x))$ for all $x \in O$. Defining $\bar{p}(s) := \sup\left\{ \frac{q(r)}{r} : 0 < r \leq s \right\}$ for $s > 0$ (well-defined since $q\left(r\right) \leq c^{-1} r$ for all $r \in \left[0, a\right]$ and since $q$ is non-decreasing which implies that $\bar{p}(s) \leq c^{-1}$ for $s \in \left(0, a\right]$ and $\bar{p}(s) = \max\left( \sup\left\{ \frac{q(r)}{r} : 0 < r \leq a \right\}, \sup\left\{ \frac{q(r)}{r} : a < r \leq s \right\} \right) \leq \max\left( c^{-1}, \frac{q(s)}{a} \right)$ for $s > a$ ), we get that $q(r) \leq r\bar{p}(r)$ for all $r > 0$. Since $\bar{p}$ is non-decreasing on $\left(0, +\infty\right)$ and non-negative we can extend $\bar{p}$ at zero by defining $\bar{p}(0) = \lim_{s \to 0^+}\left(\bar{p}(s)\right)$. Thus, we get $|x|^2 \leq W(x)\bar{p}\left(W(x)\right)$ for all $x \in O$. The proof is complete. ◁